\documentclass[11pt,leqno,draft]{amsart} 

\usepackage{mathrsfs}

\theoremstyle{plain} 

\newtheorem{global-theorem}{Theorem}
\newtheorem{theorem}{Theorem}[section]
\newtheorem{lemma}[theorem]{Lemma}

\newtheorem{corollary}[theorem]{Corollary}

\newtheorem{definition}[theorem]{Definition}
\newtheorem{proposition}[theorem]{Proposition}

\newtheorem{prop-def}[theorem]{Proposition-Definition}
\newtheorem{lemma-def}[theorem]{Lemma-Definition}

\newcommand{\eop}{\ \hfill $\Box$}

\numberwithin{equation}{subsection}

\newcommand{\uu}{{\mathbb U}}

\newcommand{\cc}{{\mathbb C}}
\newcommand{\pp}{{\mathbb P}}
\newcommand{\rr}{{\mathbb R}}

\newcommand{\zz}{{\mathbb Z}}

\newcommand{\hh}{{\mathbb H}}
\newcommand{\aaa}{{\mathbb A}}

\newcommand{\Gm}{{\mathbb G}_m}

\newcommand{\Cc}{{\mathscr C}}

\newcommand{\Dd}{{\mathscr D}}
\newcommand{\Ff}{{\mathscr F}}
\newcommand{\Oo}{{\mathscr O}}
\newcommand{\Uu}{{\mathscr U}}

\newcommand{\Mm}{{\mathscr M}}

\newcommand{\Aaa}{{\mathscr A}}

\newcommand{\Pp}{{\mathscr P}}
\newcommand{\Gg}{{\mathscr G}}
\newcommand{\Xx}{{\mathscr X}}
\newcommand{\Hh}{{\mathscr H}}
\newcommand{\Ww}{{\mathscr W}}

\newcommand{\Qq}{{\mathscr Q}}
\newcommand{\Kk}{{\mathscr K}}

\newcommand{\uQ}{{\Qq}}

\newcommand{\Fib}{{\bf Fib}}

\newcommand{\Aff}{{\bf Aff}}
\newcommand{\uHom}{\underline{\rm Hom}}

\newcommand{\Hom}{{\rm Hom}}
\newcommand{\cpx}{{\rm cpx}}
\newcommand{\wpx}{{\rm wpx}}
\newcommand{\mdl}{{\rm mod}}

\newcommand{\Perf}{{\rm Perf}}
\newcommand{\dgn}{{\bf dgn}}

\newcommand{\Spec}{{\rm Spec}}

\newcommand{\MC}{{\bf MC}}
\newcommand{\calMC}{{\Mm \Cc}}

\newcommand{\tworightarrows}{\stackrel{\displaystyle \rightarrow}{\rightarrow}}

\newcommand{\mylabel}[1]{\label{#1}}

\begin{document}

\author[C. Simpson]{Carlos Simpson}
\address{CNRS, Laboratoire J. A. Dieudonn\'e, UMR 6621
\\ Universit\'e de Nice-Sophia Antipolis\\
06108 Nice, Cedex 2, France}
\email{carlos@unice.fr}
\urladdr{http://math.unice.fr/$\sim$carlos/} 

\title[Geometricity of the Hodge filtration]{Geometricity of the Hodge filtration on the $\infty$-stack of 
perfect complexes over $X_{DR}$}

\begin{abstract}
We construct a locally geometric $\infty$-stack $\Mm _{Hod}(X,Perf)$ of perfect complexes with $\lambda$-connection structure
on a smooth projective variety $X$. 
This maps to $\aaa ^1 / \Gm $, so it can be considered as the Hodge filtration of its fiber over 1 which is $\Mm _{DR}(X,\Perf)$, 
parametrizing complexes of $\Dd_X$-modules which are $\Oo _X$-perfect. We apply the result of Toen-Vaqui\'e that $\Perf(X)$ is locally geometric. 
The proof of geometricity of the map $\Mm _{Hod}(X,\Perf) \rightarrow Perf(X)$ uses a Hochschild-like notion of weak complexes of modules over a 
sheaf of rings of differential operators. We prove a strictification result for these weak complexes, and also a strictification result for
complexes of sheaves of $\Oo$-modules over the big crystalline site. 
\end{abstract}


\maketitle



\section{Introduction}

Recall that a {\em perfect complex} is a complex of quasicoherent sheaves 
which is locally quasiisomorphic to a bounded complex of vector bundles. 
Following a suggestion of Hirschowitz, we can consider the
$(\infty ,1)$-stack $\Perf$ of perfect complexes, obtained by applying the Dold-Puppe construction to the
family of dg categories of perfect complexes considered by Bondal and Kapranov \cite{Kapranov0} \cite{BondalKapranov}. 
If $X\rightarrow S$ is a smooth projective morphism we obtain an $(\infty ,1)$-stack $\Perf (X/S)\rightarrow S$
of moduli for perfect complexes on the fibers, defined as the relative morphism stack from $X/S$ into $\Perf$.
Toen and Vaqui\'e prove
that $\Perf (X/S)$ is a {\em locally geometric} $(\infty ,1)$-stack \cite{ToenVaquie}. 
We would like to construct a locally geometric moduli stack of ``perfect complexes with integrable connection'' together with
its nonabelian Hodge filtration. 

Following Lurie \cite{LurieTopos} the notation {\em $(\infty , 1)$-category} will refer to 
any of a number of different ways of looking at $\infty$-categories whose $i$-morphisms are invertible for $i\geq 2$.
The first and easiest possibility is the notion of ``simplicial category'', however the model structure of 
\cite{Bergner} is often beneficially replaced by equivalent model category structures for Segal catgories \cite{Segal} 
\cite{Vogt} \cite{DwyerKanSmith} \cite{Tamsamani} \cite{HirschowitzSimpson} \cite{Bergner3},
quasicategories or restricted Kan complexes \cite{JoyalQuasi} \cite{BoardmanVogt}, or Rezk categories \cite{Rezk}. 
The Baez-Dolan \cite{BaezDolan}, Batanin \cite{Batanin} or any number of other $n$-category theories 
\cite{Leinster} should also apply. 
When refering to \cite{HirschowitzSimpson} for the theory of $(\infty , 1)$-stacks, we use the
model category structure for Segal $1$-categories which can be found in \cite{HirschowitzSimpson} \cite{Pellissier} \cite{Bergner3}.
We don't discuss these issues in the text. 

Throughout the paper we work over the ground field $\cc$ for convenience. The reader may often replace it by an arbitrary ground
field which should sometimes be assumed of characteristic zero. 

For a smooth projective variety $X$,
Toen has suggested that it would be
interesting to look at the notion of {\em $\Oo_X$-perfect complex of $\Dd_X$-modules}. These
can be organized into an $(\infty ,1)$-stack
in several ways which we consider and relate below. The closed points of this stack are classical objects: they are just bounded complexes
of $\Dd _X$-modules whose cohomology objects are vector bundles with integrable connection. However, the point of view of $(\infty , 1)$-stacks
allows us to consider families of such, giving a moduli object $\Mm _{DR}(X,\Perf )= \uHom (X_{DR}, \Perf )$ which generalizes the moduli stack of vector bundles
with integrable connection. Our geometricity result for this $(\infty , 1)$ moduli stack fits into the philosophy that $Hom$ stacks into geometric stacks should remain geometric under reasonable
hypotheses \cite{geometricN} \cite{Aoki}. 

One motivation for considering perfect complexes is that they arise naturally in the context of variation of cohomology
studied by Green-Lazarsfeld \cite{GreenLazarsfeld}. This corresponds exactly to our situation. 
Given a smooth morphism $f: X\rightarrow Y$, even if we start 
with a family of vector bundles with integrable connection $\{ E_s\} _{s\in S}$ on $X$, the family of  higher direct images $\rr f_{\ast}(E_s)$
will in general be a perfect complex of $\Dd _X$-modules on $Y\times S/S$. 
The Gauss-Manin connection on the full higher derived direct
image was considered by Katz and Oda \cite{Katz} \cite{KatzOda}, see also Saito \cite{Saito89}. 
The higher direct image may profitably be seen as a morphism $\rr f_{\ast} : \Mm _{DR}(X,\Perf )\rightarrow \Mm _{DR}(Y,\Perf )$
between $(\infty ,1)$ moduli stacks. 

Following a suggestion of Deligne, the de Rham moduli spaces deform into corresponding Dolbeault moduli spaces of Higgs bundles
\cite{naht} \cite{hodgefilt}.
This deformation can be viewed as the {\em nonabelian Hodge filtration} generalizing \cite{NeisendorferTaylor}, \cite{Deninger} 
in the abelian case. For moduli of perfect complexes, we get a
stack $\Mm _{Hod}(X,\Perf )\rightarrow \aaa ^1$ with $\Gm$-action, whose fiber over $\lambda =1$ is $\Mm _{DR}(X,\Perf )$
and whose fiber over $\lambda =0$ is the moduli stack $\Mm _{Dol}(X,\Perf )$ of $\Oo$-perfect complexes of Higgs sheaves whose cohomology sheaves
are locally free, semistable, with vanishing Chern classes. The family $\Mm _{Hod}(X,\Perf )$ is the Hodge filtration on the moduli of perfect
complexes over $X_{DR}$.  

Deligne's original idea of looking at the moduli space of vector bundles with $\lambda$-connection was designed to provide a construction of
Hitchin's twistor space \cite{Hitchin}. We can imagine making the same glueing construction with $\Mm _{Hod}(X,\Perf )$ to construct a twistor space of
perfect complexes. Investigation of this kind of object fits into a long running research project with 
L. Katzarkov, T. Pantev, B. Toen, and more recently, G. Vezzosi,
M. Vaqui\'e, P. Eyssidieux, about nonabelian mixed Hodge theory \cite{namhs}. The moduli of perfect complexes provides an interesting special
case in which the objects being considered are essentially linear. Even in the nonlinear case perfect complexes will appear as linearizations,
so it seems important to understand this case well. 

Our purpose in the present paper is to show that the objects entering into the above story have a reasonable geometric structure.
The notion of {\em geometric $n$-stack} is a direct generalization of Artin's notion of algebraic stack \cite{Artin}, and a locally geometric $\infty$-stack
is one which is covered by geometric $n$-stacks for various $n$. In our linear case, the objects are provided with morphism complexes containing
morphisms which are not necessarily invertible, thus it is better to speak of $(\infty ,1)$-stacks. 
\newline
---Make the convention \label{conven} that the word {\em stack} means {\em $(\infty ,1)$-stack} when that makes sense.

Our main result is that the moduli stack
$\Mm _{Hod}(X,\Perf))$ is locally geometric \ref{thethm} \ref{hfthm}. Its  proof relies heavily on the recent result of Toen and Vaqui\'e
that $\Perf(X)$ is locally geometric. We are thus reduced to proving 
that the morphism
$$
\Mm _{Hod}(X,\Perf) \rightarrow \Perf (X\times \aaa ^1 / \aaa ^1) = \Perf (X)\times \aaa^1
$$
is geometric.

It seems likely that geometricity could be deduced from Lurie's representability theorem, and might also
be a direct consequence of the formalism of Toen-Vaqui\'e. Nonetheless, it seems interesting to have a reasonably
explicit description of the fibers of the map: this means that we fix a perfect complex of $\Oo$ modules $E$ over $X$
and then describe the possible structures of $\lambda$-connection on $E$. 
The notion of $\lambda$-connection is encoded
in the action of a sheaf of rings of differential operators $\Lambda_{Hod}$, which gives back $\Dd_X$ over $\lambda = 1$. 

If $f:X\rightarrow Y$ is a smooth morphism, 
the higher direct image gives a morphism of locally geometric stacks
$$
\rr f_{\ast} : \Mm _{Hod}(X,\Perf)\rightarrow \Mm _{Hod}(Y,\Perf),
$$
which is a way of saying that the higher direct image functor between de Rham moduli stacks preserves the Hodge filtration.

We work in the more general situation of a {\em formal category of smooth type} $(X,\Ff )$ which was introduced in 
crystalline cohomology theory \cite{dixexposes} \cite{Berthelot} \cite[Ch. VIII]{Illusie} as a generalization of de Rham theory. 
Associated to $(X,\Ff )$ is a sheaf of rings of differential operators $\Lambda$. 
The classical case is when $X$ is smooth and $\Ff _{DR}$ is the formal completion of the diagonal in $X\times X$.
Then $\Lambda _{DR}=\Dd _X$. 
Deformation to the normal cone along the diagonal provides a deformation $(X\times \aaa ^1, \Ff _{Hod})\rightarrow \aaa ^1$
embodying the Hodge filtration. 
The corresponding ring $\Lambda _{Hod}$
is the Rees algebra of $\Dd _X$ for the filtration by order of differential operators. 

Working in the context of a formal category has the advantage of streamlining notation, and potentially could give applications
to other situations such as Esnault's $\tau$-connections along foliations \cite{Esnault}, logarithmic connections 
\cite{DeligneRS} \cite{Burgos} \cite{Nitsure} \cite{NitsureSabbah} and so forth. 
Some interesting properties of perfect complexes with logarithmic connections were suggested by conversations with Iyer, see the appendix to the preprint
version of \cite{IyerSimpson}; this was one of the motivations for the present work. 

Our description of the space of $\Lambda$-module structures on $E$ passes through a Kontsevich style  
Hochschild weakening of the
notion of complex of $\Lambda$-modules. In brief, the tensor algebra 
$$
T\Lambda := \bigoplus \Lambda \otimes _{\Oo_X} \ldots \otimes _{\Oo_X}\Lambda
$$
has a differential and coproduct, and for $\Oo$-perfect complexes $E$ and $F$, this allows us to define the complex 
$$
Q(E,F):= Hom (T\Lambda \otimes _{\Oo _X}E,F)
$$
with composition. A weak structure is an element $\eta \in Q^1(E,E)$ satisfying the {\em  Maurer-Cartan equation} $d(\eta )+\eta ^2=0$.
The basic technique for working with these objects comes from the paper of Goldman and Millson \cite{GoldmanMillson}.
This procedure works on affine open sets, and we need a \v{C}ech globalization---again using Maurer-Cartan---to get to a projective $X$. 

The idea that we have to go to weak structures \cite{Stasheff} \cite{GugenheimLambe} in order to obtain a good
computation, was exploited by Kontsevich \cite{Kontsevich} \cite{BarannikovKontsevich}
\cite{KontsevichSoibelman}, see also Tamarkin-Tsygan \cite{TamarkinTsygan},  Dolgushev \cite{Dolgushev}, Borisov \cite{Borisov} and others. 
Kontsevich used this idea to calculate the deformations of various objects by
going to their weak versions (e.g. $A_{\infty}$-categories). Hinich also cites Drinfeld, and historically
it also goes back to Stasheff's $A_{\infty}$-algebras, Toledo and Tong's twisted complexes, the theory of operads and many other
things. 
The tensor algebra occurs in connection with such weak structures in Johansson-Lambe \cite{JohanssonLambe}. 
The combination of a Hochschild-style differential and a coproduct occurs in Pridham's notion of {\em simplicial deformation complex}
\cite{Pridham} and much of what we do here follows pretty much the same kind of formalism. 
Looking at things in this way was
suggested to me by E. Getzler, who was describing his way of looking at some other related questions \cite{Getzler}, and also in a conversation about
small models with V. Navarro-Aznar. 

The use of the Hochschild complex for encoding weak de Rham structures played a central role in the construction of characteristic
classes by Block and Getzler \cite{Block} \cite{Getzler} \cite{BlockGetzler}, Bressler, Nest, and Tsygan \cite{BresslerNestTsygan}   \cite{NestTsygan1} \cite{NestTsygan3}, 
Dolgushev \cite{Dolgushev}. This goes back to some of the old sources of homological deformation theory \cite{Gerstenhaber}.  
It has been related to de Rham cohomology and the Hodge to de Rham spectral 
sequence by Kontsevich, Kaledin \cite{Kaledin}, and to the theory of loop spaces by Ben-Zvi and Nadler
\cite{BenZviNadler}. The Rees algebra construction plays an important role in these theories. 
These things are related to ``deformation quantization'' \cite{KontsevichDQ} \cite{KashiwaraSchapira}.  

The application
to weak $\Lambda$-module structures is a particularly easy case since everything is almost linear (i.e. there are no
higher product structures involved). Our argument is structurally similar to Block-Getzler \cite{BlockGetzler}. 
An important step in the argument is the calculation of the homotopy fiber product involved in 
the definition of geometricity. This is 
made possible by the combination of Bergner's model category structure on the category of simplicial categories \cite{Bergner}
and Tabuada's model structure for dg-categories \cite{Tabuada}.  

Another route towards parametrizing complexes of $\Dd _X$-modules is considered extensively by Beilinson and
Drinfeld  \cite{BeilinsonDrinfeld} \cite{BeilinsonDrinfeld2}, based on the description as modules over the dg algebra
$(\Omega ^{\cdot}_X, d)$ by
Herrera and Lieberman \cite{LiebermanHerrera}, see also Kapranov \cite{Kapranov2} and Saito \cite{Saito89}. It should be possible to prove the 
geometricity theorem
using this in place of the weak $\Lambda$-module description, and these notions enter into \cite{BenZviNadler}. 

Toen' suggestion that perfect complexes over $X_{DR}$ remain interesting objects,
essentially because they encode higher cohomological data, fits in with his notion of
{\em complex homotopy type} $X\otimes \cc$. The basic idea is that tensor product provides the 
$\infty$-category of perfect complexes over $X_{DR}$ with a Tannakian structure, and $(X\otimes \cc )_{DR}$
is the Tannaka dual of this Tannakian $\infty$-category.  It continues and generalizes rational homotopy theory 
\cite{DGMS} \cite{Tanre} \cite{Gomez-Tato} \cite{Hain2}. 
This  Tannakian theory could be seen as another motivation for looking at the locally geometric stacks 
$\Mm _{DR}(X,\Perf )$ we consider here: tensor product provides a monoidal structure on $\Mm _{DR}(X,\Perf)$ and a basepoint
$x\in X$ gives a fiber-functor 
$\omega _{x,DR} : \Mm _{DR}(X,\Perf) \rightarrow \Perf $.
It would be interesting to  investigate further this structure in view of the geometric property of $\Mm _{DR}(X,\Perf)$. 

Katzarkov, Pantev and Toen provide the schematic homotopy type $X\otimes \cc$
with a mixed Hodge structure \cite{KaPaTo1} \cite{KaPaTo2}. This leads to Hodge-theoretic
restrictions on the homotopy type of $X$ generalizing those of 
\cite{DGMS}, and should be related to the nonabelian mixed Hodge structure
on the formal completions of $\Mm _{DR}(X,\Perf )$ in much the same way as Hain's mixed Hodge structure on
the relative Malcev completion \cite{Hain2} is related to the mixed Hodge structure on the formal completion of
$\Mm _{DR}(X,GL(n))$. 
Conjecturally, the local structure of the Deligne-Hitchin twistor space associated to $\Mm _{DR}(X,\Perf )$
should be governed by the theory of
mixed Hodge modules \cite{SaitoMHM} \cite{MSaito0504} or mixed twistor modules \cite{Sabbah}.

The notion of connection on higher homotopical structures can be viewed as an example of {\em differential geometry of gerbes}, a theory
extensively developped by Breen and Messing \cite{BreenMessing1} \cite{BreenMessing2}, related to many other things such
as \cite{Navarro-AznarGM} \cite{Picken} \cite{BlockGetzler} \cite{Block05} \cite{Lott} \cite{AschieriCantiniLurco} \cite{Behrend}
\cite{SchreiberStrobl}. In this point of view, higher categorical objects over a formal category $(X,\Ff )$ can be viewed as given by ``infinitesimal cocycles'' \cite{AKock} 
where the transition functions correspond to the infinitesimal morphisms in $\Ff$. While we don't explicitly mention this idea in our argument, it 
represents the underlying intuition and motivation.

We look only at sheaves of rings of differential operators which come from formal categories of smooth type.
This condition insures the almost-polynomial hypothesis crucial for bounding the $Ext$ groups via the
HKR theorem \cite{HKR} \cite{Gerstenhaber} \cite{Dolgushev}, 
and also insures the existence of bounded de Rham-Koszul resolutions. 
It would undoubtedly be interesting to treat various more general sheaves of rings $\Lambda$.
That might be useful, for example in the case of 
Azumaya algebras \cite{Lieblich} \cite{HoffmanStuhler}.
Furthermore, it would probably be a good idea to go all the way towards looking at complexes of
modules over a sheaf of DG-algebras \cite{Hinich} \cite{Block05} or even a DG-stack \cite{Tabuada} \cite{Porta}, 
and that might help with implementing the proof using \cite{LiebermanHerrera} \cite{BeilinsonDrinfeld} but we don't do that here.

We prove two main strictification results. These proofs were absent from the first preprint version and constitute one of the
major changes here. The first result concerns the passage between $\Lambda$-modules and sheaves on the big crystalline site
$CRIS(X,\Ff )$. A complex of $\Lambda$-modules $E$ leads to a complex of sheaves of $\Oo$-modules $\Diamond (E)$
on $CRIS(X,\Ff )$. The image of this functor consists of the {\em quasicoherent} complexes
on $CRIS(X,\Ff )$, the complexes whose pullback morphisms are isomorphisms.  However, a general complex of sheaves of $\Oo$-modules 
is very far from satisfying this property. We show that $\Oo$-perfect complexes on $CRIS(X,\Ff )$ are quasiisomorphic to quasicoherent
ones, that is to say that $\Diamond$ is an equivalence between the $(\infty ,1)$-categories of $\Oo$-perfect objects on both sides
\ref{coherence1}, \ref{diamondequiv}. 
The proof over affine $X$ is a standard argument using the adjoint $\Diamond ^{\ast}$ inspired by some
remarks in \cite{BerthelotOgus}. 

Before that, we also show that the $(\infty ,1)$-category of $\Oo$-perfect complexes on $CRIS(X,\Ff )$
is equivalent to the $(\infty ,1)$-functor category from the stack $[X/\! /\Ff ]$ corresponding to our formal category, to
$\Perf$. I would like to thank Breen for pointing out that this passage between nonabelian cohomology of $[X/\! /\Ff _{DR}]$
and usual de Rham theory as manifested in $\Dd _X$-modules, wasn't entirely obvious and needed careful consideration. 
A cocycle version of this kind of statement is provided by the works of Breen and Messing \cite{BreenMessing1} \cite{BreenMessing2}. 

The second main strictification result says that a Hochschild-style weak complex of $L$-modules $(E,\eta )$, 
determined by a Maurer-Cartan element as described above,
comes up to homotopy from a complex of strict $L$-modules. 
The weak  notion has a key homotopy invariance property as in 
Gugenheim-Lambe-Stasheff-Johansson \cite{GugenheimLambe} \cite{JohanssonLambe}: if $E$ is replaced by an equivalent complex of $A$-modules $E'$,
then $\eta$ can be modified to an equivalent MC element $\eta '$ for $E'$ (Lemma \ref{maininvariance}). This property exactly says 
(Corollary \ref{Pfibration}) that
the functor from the dg category of weak complexes, to the dg category of perfect complexes over $A$, is fibrant in the sense of
Tabuada \cite{Tabuada}, hence its Dold-Puppe is fibrant in the sense of Bergner \cite{Bergner}. This is the technique which allows us to 
compute the homotopy fiber of the map $\Mm (\Lambda ,\Perf)\rightarrow \Perf (X)$ (see also \cite{Tabuada2}), 
to go from the Toen-Vaqui\'e geometricity result \cite{ToenVaquie}
for $\Perf (X)$ to geometricity of $\Mm (\Lambda ,\Perf )$. 

Both strictification results are first proven over affine open sets, and indeed the functors in question are first defined in
the affine case. We then pass to the global $X$ by using the notion of stack (that is, $(\infty ,1)$-stack or Segal $1$-stack) \cite{HirschowitzSimpson}, 
or a similar dg version of descent.
This descent technique goes back to Toledo-Tong's twisted complexes \cite{ToledoTong}, and fundamentally speaking comes from
the cohomological descent technique pioneered by Deligne and Saint-Donat. 

One way to see why we need to globalize by descent is to consider the problem of defining appropriate model category
structures. The most widely usable model category structure for complexes is Hovey's injective model structure \cite[arxiv version]{Hovey} on the
category of unbounded complexes $\cpx (\Aaa )$ in any AB5 abelian category $\Aaa$. However, this structure is not appropriate for
left-derived functors. For example the crystallization functor $\Diamond$ takes a complex of $\Lambda$-modules $E$ to a sheaf
on the big site $CRIS(X,\Ff )$ whose value over $\Spec (B)$ is $E\otimes _{\Oo (X)}B$. The presence of this tensor product
means that $\Diamond$ is a left derived functor, so we would want it to start from a projective  model category structure.
However, the category of $\Lambda$-modules on a global $X$ doesn't have enough projectives. 

This is a well-known problem in the
theory of $\Dd _X$-modules which I first learned about in J. Bernstein's course. This kind of problem has been 
attacked by Deligne at various times, leading to notions of pro and ind objects in various categories. Fiorot's recent
paper \cite{Fiorot} gives a thorough explanation of Deligne's theory of stratified or crystalline pro-modules. It is likely
that this direction could lead to a more natural description of the global objects which we have considered above.
We meet a similar kind of structure in \S \ref{finiterep} at the end of our proof. 
 
Deligne's original cohomological descent techniques have led to many developments such as in \cite{GuillenNavarro}
\cite{Wojtkowiak}. 
It is worth
noting that the basic construction for descent in \cite[19.4]{HirschowitzSimpson} is the same as the simplicial descent construction used in 
SGA IV and subsequently \cite{hodge3}. 
\v{C}ech descent is closely related to Hinich's work \cite{Hinich} \cite{Hinich2} as well as to the twisted complexes
of Toledo and Tong \cite{ToledoTong}, see also Siegel \cite{Siegel}, Beke \cite{Beke}, Block \cite{Block05}. 
Stasheff and Wirth recently posted an article based on Wirth's 1965 thesis
\cite{StasheffWirth} concerning the cocycle description of descent data for higher stacks; this kind of cocycle
condition plays an important part in \cite{Breen} \cite{BreenMessing1} \cite{BreenMessing2}.
These are all examples of {\em homotopy-coherent category theory} \cite{CordierPorter} \cite{SchwanzlVogt}.
Several authors have considered this type of descent question for Deligne cohomology \cite{DupontLjungmann} \cite{SchreiberStrobl}.
Another very recent
posting concerning the \v{C}ech globalization of Maurer-Cartan elements for a sheaf of DGLA's is \cite{FMM}. 
Tabuada has explained to me that the \v{C}ech-style expression for a limit of dg-categories coming out of 
\cite{ToledoTong} is also closely related to an expression for the homotopy path categories for functors of
dg categories in \cite{BondalKapranov} \cite{Tabuada}. 

Ideally the dg \v{C}ech globalization should fit into a theory of {\em dg stacks} based on the model category structure of \cite{Tabuada},
having \cite{ToledoTong} as the main example,
and compatible with the simplicial notion via Dold-Puppe. 
A similar remark which could have been made earlier is that 
the \v{C}ech globalization of a presheaf of d.g.a.'s is the homotopy limit in Hinich's model category \cite{Hinich}.
In that context
the $\widetilde{DP}\MC$  construction \S \ref{mcstack} should be viewed as basically the same as the simplicial Deligne groupoid constructed by Hinich \cite{Hinich1}. 
In the Hochschild-like context  of weak complexes, it is also closely related to
the work of Pridham on ``simplicial deformation complexes'' \cite{Pridham}.

The original preprint version of this paper \verb+arXiv:math/0510269v1+ contains some extra introductory motivation concerning
$\lambda$-connections and the Hodge filtration,
some of it corresponding more closely to my talk at the Deligne conference in Princeton, October 2005. 
The present version aims to provide a fuller account of the proofs. The proof of the geometricity theorem follows
the same lines except that the trick of using the Rees direct sum construction in the original preprint has been replaced here
by a more careful consideration of the filtration. Newly added here is a discussion of the passage between complexes of $\Lambda$-modules
and sheaves on the crystalline site associated to a formal category \S \ref{strictif-crys}. This passage can be skipped if the reader wants to interpret
the geometricity result as speaking directly about the moduli stack of complexes of $\Lambda$-modules. 

One of our main long range goals is to develop and generalize mixed Hodge theory to a nonabelian setting, guided by Deligne's outlook.
He has led us to understand Hodge theory as one of the most visible manifestations of the special topological structure of
algebraic varieties, going alongside other manifestations such as the action of Frobenius on $\ell$-adic cohomology. This point of view imposes that we consider Hodge theory for any kind of nonabelian cohomology or indeed any other topological structure, so it came
as no surprise that Deligne's interpretation of Hitchin's quaternionic structure led immediately to the notion of nonabelian Hodge filtration. It is much less immediate to understand the role of the weight filtration in this picture, although there is no doubt that it is there. If we haven't considered this question directly in the present paper, the hope is that the present example of
moduli of perfect complexes will prove to be an interesting one from the point of view of mixed Hodge theory. 
On a technical level, much of the inspiration for this work comes from the ideas which Deligne has transmitted, to others as well as myself, in the form of letters.
A great many of our numerous references themselves reflect the influence of his ideas. So, I would like to thank
Deligne for showing the way with generosity and encouragement. Then I would also like to thank a number of my co-inspirees 
for helpful conversations and remarks some of which are mentionned in the text.


\section{Differential graded model categories of complexes}

Since we are working with complexes, the notion of dg enrichment \cite{Kapranov0} \cite{BondalKapranov} provides
a good approach to the notion of $(\infty ,1)$ category. The relationship between complexes and simplicial sets \cite{Illusie}
is given by Dold-Puppe.

\subsection{The Dold-Puppe construction}

Define the following complex denoted $D(n)$. Let 
$D(n)^{-k}$ denote the free abelian group generated by the inclusions 
$\varphi :[k] \hookrightarrow [n]$, for $0\leq k \leq n$, with $D(n)^{-k}=0$ otherwise. 
Denote the generator corresponding to $\varphi$ by just $\varphi$.

Define the differential $d:D(n)^{-k} \rightarrow D(n)^{1-k}$
by 
$$
d\varphi := \sum _{i=0}^k(-1)^i (\varphi \circ \partial _i)
$$
where $\partial _i : [k-1]\rightarrow [k]$ is the face map skipping the $i$th object. 

Define a coproduct $\kappa : D(n)\rightarrow D(n)\otimes D(n)$ as follows. 
Let $lt_i : [i]\rightarrow [k]$ denote the inclusion of the first $i+1$ objects,
and $gt_i : [k-i]\rightarrow [k]$ the inclusion of the last $k+1-i$ objects. 
Note that these overlap, both including object number $i$.
For $\varphi : [k]\hookrightarrow [n]$, put
$$
\kappa (\varphi ) := \sum _{i=0}^{k} (\varphi \circ lt_i) \otimes (\varphi \circ gt_i).
$$
This is co-associative and  compatible with the differentials.

Our collection of complexes $D(n)$ forms a cosimplicial object in the category of complexes.
Namely, for any map $\psi : [n]\rightarrow [m]$ we obtain a map $D(n)\rightarrow D(m)$ sending 
$\varphi $ to $\psi \circ \varphi$ if the latter is injective, and to zero otherwise. 

We can now define the {\em Dold-Puppe functor} which to a complex $A$ associates the simplicial 
abelian group which in degree $n$ has
$$
DP(A)_n:= \Hom (D(n),A)
$$
(the group of morphisms of complexes from $D(n)$ to $A$). This is a simplicial object by functoriality of
$D(n)$ in $n$. Elements of $DP(A)_n$ may be considered as functions $a(\varphi )\in A_k$ for each
$\varphi : [k]\hookrightarrow [n]$, satisfying $d(a(\varphi )) = \sum (-1)^ia(\varphi \circ \partial _i)$.

On a theoretical level, let $N$ denote the normalized complex functor from simplicial abelian groups to
nonpositively graded complexes of abelian groups. If we let $\dgn$ denote the degenerate subcomplex then $N(A)={\bf s}(A)/\dgn (A)$. 
The construction $DP$ restricted to negatively graded complexes is the inverse of $N$. On unbounded complexes
$DP=DP\circ \tau _{\leq 0}$ is right adjoint to $N$ or more precisely to $\tau _{\leq 0}^{\ast}N$ where the left adjoint of the intelligent truncation $\tau _{\leq 0}$
is just inclusion of negatively graded complexes into unbounded complexes.   

Compose $\Hom (D(n),A)\times \Hom (D(n),B) \rightarrow \Hom (D(n)\otimes D(n),A\otimes B)$
with the the map of right composition by the coproduct $\kappa$ to get 
$$
\Hom (D(n),A)\times \Hom (D(n),B) \rightarrow \Hom (D(n),A\otimes B).
$$
We obtain maps of simplicial sets (which however are not linear)
$$
DP(A)\times DP(B)\rightarrow DP(A\otimes B).
$$
This product is associative, because of co-associativity of $\kappa$. Note also that if $\zz [0]$ denotes the 
complex with $\zz$ in degree $0$ then $DP(\zz [0])$ is the constant simplicial abelian group $\zz$. 
With respect to this isomorphism, the product above is unital too.

If $R$ is a commutative base ring and $A$ and $B$ are complexes of $R$-modules then we have a map
$A\otimes B\rightarrow A\otimes _RB$. The product
$$
DP(A)\times DP(B)\rightarrow DP(A\otimes _RB),
$$
remains associative and unital.

These facts allow us to use $DP$ to create a simplicial category out of a differential graded category.
A differential graded category over a ring $R$ is a category $C$ enriched in complexes of $R$-modules.
Applying the functor $DP$ to the enrichment, with the product maps constructed above, yields 
a category which we denote by $\widetilde{DP}(C)$. The objects of this category are the same as the objects of
$C$, whereas if $x,y\in {\rm ob}(C)$ then by definition 
$$
Hom _{\widetilde{DP}(C)}(x,y):= DP(Hom_C(x,y)).
$$
The composition of morphisms is given by the above product maps:
$$
\Hom^{\cdot}_C(x,y)\otimes _R \Hom^{\cdot}_C(y,z)\rightarrow \Hom^{\cdot}_C(x,z)
$$
yields
$$
DP(\Hom^{\cdot}_C(x,y))\times DP(\Hom^{\cdot}_C(y,z))\rightarrow 
DP(\Hom^{\cdot}_C(x,z)).
$$
The 
morphisms of $\widetilde{DP}(C)$, i.e. the points in the degree $0$ part of the simplicial mapping sets,
are the same as the morphisms of $C$. A morphism of $C$ is an inner equivalence if and only if its image
is an equivalence in $\widetilde{DP}(C)$.

\subsection{Differential-graded enrichment for model categories}
\mylabel{sec-dg-enrich}

Quillen formulated a notion of ``simplicial model category'' \cite{Quillen}. In some important examples, this provided an
easy way to obtain a good simplicial category. Later the construction of a simplicial category in general was
provided by the notion of Dwyer-Kan localization \cite{DwyerKan}. For the problem of differential-graded enrichment, Drinfeld introduced 
a general localization procedure \cite{Drinfeld}. However, as with simplicial enrichment, many basic examples 
can be treated by an easier route generalizing Quillen's 
definition to a notion of {\em dg model category}.

Suppose $\Cc$ is a closed model category which also has a structure of $\cc$-linear dg-category.
For $P\in {\rm cpx}(\cc )$, an {\em object $X\otimes P$} is an object of $\Cc$ together with a map $P\rightarrow Hom ^{\cdot}(X,X\otimes P)$
universal for maps $P\rightarrow Hom ^{\cdot}(X,Y)$. Dually an {\em object $Y^P$} as an object of $\Cc$ together
with a map $P\rightarrow Hom^{\cdot}(Y^P,Y)$ universal for maps $P\rightarrow Hom ^{\cdot}(X,Y)$. 

Adapt the numbering of \cite{Quillen} to the dg case. 
Say that the condition $(DGM0)$ holds if $X\otimes P$ and $Y^P$ exist for every complex $P$ and objects $X,Y\in \Cc$.
If this is the case, then the
associations $X,P\mapsto X\otimes P$ and $X,P\mapsto X$ are functorial in $P$ and $X$. 
They are adjoints to the $Hom$ complex functors, and the $Hom$ complex functor is  
compatible with colimits and limits
in $\Cc$ in the sense that $X,Y\mapsto Hom^{\cdot}(X,Y)$ takes colimits in the variable $X$ to limits of complexes,
and takes limits in the variable $Y$ to limits of complexes.

Say that $\Cc$ together with its enrichment satisfies condition $(DGM7)$ if for any cofibration $i:U\rightarrow V$
and cofibration $p:X\rightarrow Y$, the map 
$$
Hom^{\cdot}(V,X)\rightarrow Hom^{\cdot}(U,X)\times _{Hom^{\cdot}(U,Y)}Hom^{\cdot}(V,Y)
$$
is a fibration of complexes (i.e. a surjection) which is trivial if either $i$ or $p$ is trivial. We say that the closed model category
$\Cc$ together with its dg-enrichment is a {\em dg model category} (over $\cc$) if it satisfies conditions $(DGM0)$ and 
$(DGM7)$. 

I would like to thank Toen for pointing out that the
definition above is different from that used in 
\cite{HoveyBook} and \cite{ToenMorita}. Here there is a natural map 
$$
{\bf a}:X\otimes (P\otimes Q)\rightarrow (X\otimes P)\otimes Q
$$
satisfying an associativity condition, however this natural map is not necessarily an isomorphism
as required in the definition of \cite{HoveyBook}.

It will be convenient to go in the other direction and use a construction of the form 
$X,P\mapsto X\otimes P$ to create the dg enrichment. 
A {\em lax dg module structure} on a $\cc$-linear model category $\Cc$
is a functorial product $U,P\mapsto U\otimes _kP$ for $U\in \Cc$ and $P\in {\rm cpx}(k)$, together with
a natural map ${\bf a}$ as above, not necessarily an isomorphism but satisfying an associativity axiom. 

Given this product, we can define a differential graded enrichment of $\Cc$. For $U,V\in \Cc$ let
$\Hom^{\cdot} (U,V)$ be the complex of $k$-vector spaces defined by adjunction: a map
$$
P\rightarrow \Hom^{\cdot} (U,V)
$$
of $k$-complexes is the same thing as a map $U\otimes P\rightarrow V$. If we let $P_n$ be the complex $k\rightarrow k$
with the first term in degree $n$, then a map $P_n\rightarrow \Hom^{\cdot} (U,V)$ is the same as an element of degree $n$
in $\Hom^{\cdot} (U,V)$, thus we can {\em define}
$$
\Hom^n(U,V) := {\rm Mor}_{\Cc}(U\otimes P_n, V).
$$
The morphism $P_{n+1}\rightarrow P_n$ yields the differential on the complex $\Hom^{\cdot} (U,V)$ whose square is zero. 

We require that the product $U,P\mapsto U\otimes P$ preserve colimits in ${\rm cpx}(k)$ in the second variable. This insures
that the complex defined above is adjoint to the tensor product. In this case we say that the dg module structure induces
a dg enrichment, with composition maps determined by ${\bf a}$.

Suppose the functor $U,P\mapsto U\otimes P$ preserves cofibrations and trivial cofibrations in both variables. 
Then the dg enriched structure on $\Cc$ defined above makes it into a dg model category, with tensor product operation
being the same as the given one.

The following proposition is a model category version of the construction of Kapranov \cite{Kapranov0} and Bondal-Kapranov \cite{BondalKapranov}.

\begin{proposition}
\mylabel{dgsimpl}
If $\Cc$ is a dg model category, then $\widetilde{DP}\Cc$ is a simplicial model category with the same underlying category,
the same weak equivalences, cofibrations and fibrations. If $\Cc _{\rm cf}$ is the dg subcategory of fibrant and cofibrant objects
then we can calculate the Dwyer-Kan simplicial localization as $L\Cc \cong \widetilde{DP}\Cc_{\rm cf}$.
\end{proposition}
\begin{proof}
If $K$ is a simplicial set and $U\in \Cc$ put
$U\otimes K:= U\otimes N(\cc [K])$. The multiplicativity property of $DP$ is adjoint to 
$$
N(A\otimes B) \rightarrow N(A)\otimes N(B)
$$
and this transforms the lax associativity constraint ${\bf a}$ for $\Cc$ into a constraint defining a lax simplicial module structure.
The associated simplicial category is $\widetilde{DP}\Cc$. The functor $K\mapsto N(\cc [K])$ preserves cofibrations and trivial cofibrations
so we get the required property $SM7(b)$. For the calculation of the localization note that in a Quillen simplicial model category,
$L\Cc$ is equivalent to the simplicial
category  of cofibrant and fibrant objects \cite{DwyerKan}. 
\end{proof}

\subsection{Homotopy fibers}
\mylabel{htyfib}

Tabuada constructs a good model category of differential graded categories \cite{Tabuada}.
This allows us to consider the homotopy fiber of a functor of d.g.c.'s, and via Dold-Puppe it is compatible with
taking the homotopy fiber of a functor of simplicial categories. Indeed, the fibrancy condition for d.g.c.'s \cite{Tabuada}
\cite{ToenMorita}
corresponds nicely with the Dwyer-Kan-Bergner fibrancy condition for simplicial categories \cite{Bergner}.
Once we have a fibrant morphism, the homotopy fiber is the same as the usual strict fiber.
Recently Tabuada has given a separate proof of this relationship in 
\cite{Tabuada2}. 

In Tabuada's model structure \cite{Tabuada} \cite{ToenMorita} a dg functor 
$f:A\rightarrow B$ is {\em fibrant} if it satisfies the following two conditions:
\newline
(1)\,\, for each pair of objects $x,y$ of $A$, the map 
$$
Hom _A(x,y) \rightarrow Hom _B(f(x), f(y))
$$
is a surjection of complexes of abelian groups; and
\newline
(2) \,\, for $x\in {\rm ob}(A)$ and $y\in {\rm ob}(B)$ and $u$ an equivalence in $B$ 
between $f(x)$ and $y$, then $u$ lifts to an equivalence
$\tilde{u}$ in $A$, between $x$ and a lift $\tilde{y}$ of $y$.

If $f$  is fibrant in the above sense, then $\widetilde{DP}(f)$ is a functor between simplicial categories
which is fibrant in Bergner's model structure \cite{Bergner}, namely
\newline
(1) \,\, $\widetilde{DP}(f)$  induces a Kan fibration on mapping complexes; and
\newline
(2) \,\, equivalences in $\widetilde{DP}(B)$ lift to $\widetilde{DP}(A)$ with one endpoint fixed in the same
way as above. 

In particular, once we have Bergner's fibrancy condition,
we can use $\widetilde{DP}(f)$ directly to form the homotopy fiber product:
$$
\Fib (\widetilde{DP}(f) / b) := \widetilde{DP}(A)\times _{\widetilde{DP}(B)} \{ b\} .
$$ 
The objects of $\Fib (\widetilde{DP}(f) / b)$ are the objects of $A$ which map to $b$,
and the simplicial mapping sets are the subsets of the mapping sets in $\widetilde{DP}(A)$
of objects mapping to (the degeneracies of) the identity of $b$. 

We can construct the fiber on the level of differential graded categories. However, the fiber is not itself
a Dold-Puppe construction but only a closely related affine modification, because of the condition that
the morphisms map to the identity of $b$. In order to define this structure, first say that an 
{\em augmented differential graded category} $(C,\varepsilon )$ is a dgc $C$ with maps
$$
\varepsilon : Hom ^0_C(x,y)\rightarrow \cc
$$
for all pairs of objects $x,y$, such that $\varepsilon (dv)= 0$ when $v$ has degree $-1$, and 
such that $\varepsilon (uv)= \varepsilon (u) \varepsilon (v)$ and $\varepsilon (1_x)=1$. 
It is the same thing as a functor to the
dgc $(\ast , \cc [0])$ of one object whose endomorphism algebra is $\cc$ in degree $0$. 

If $(A,\varepsilon )$ is an augmented dgc, define the {\em affine Dold-Puppe} 
\label{affineDP}
$\widetilde{DP}(A,\varepsilon )$ to be the
sub-simplicial category of $\widetilde{DP}(A)$ whose simplices are those which project to degeneracies
of the unit $1\in \cc$ in $DP(\cc [0])$, this latter being the constant simplicial group $\cc _{\Delta}$. 

If $f$ is fibrant in the above sense, define an augmented dgc $\Fib ^{+}(f/b)$ as follows.
The objects are the objects of $A$ which map to $b$. The mapping spaces are the subcomplexes of
$Hom _A(x,y)$ consisting of elements which map to $0$ in degree $\neq 0$ and which map to a constant multiple
of the identity $1_b$ in degree $0$. The augmentation $\varepsilon$ is the  map to the complex line of multiples of 
$1_b$. 

\begin{lemma}
\label{dpfiber}
If $f:A\rightarrow B$ is a functor of differential graded categories which is fibrant in the above sense,
and if $b\in {\rm ob}(B)$,
then the homotopy fiber of $\widetilde{DP}(f)$ over $b$ is calculated by the affine Dold-Puppe
of the augmented dgc fiber of $f$,
$$
\Fib (\widetilde{DP}(f)/ b) = \widetilde{DP}(\Fib ^{+}(f/b),\varepsilon ).
$$
\end{lemma}
\eop

\subsection{Injective and projective model structures on categories of complexes}
\mylabel{sec-catcomplex}

It is now possible to work with unbounded complexes, using model category theory
\cite{AJS} \cite{Hovey} \cite{Hovey} \cite{Serpe} \cite{Spaltenstein} \cite{JardineCpx}.
This is not strictly speaking necessary for what we do, but it allows us to avoid cumbersome 
consideration of upper and lower boundedness conditions. 

Let $\Aaa$ denote an AB5 abelian category. 
Let ${\rm cpx}(\Aaa )$ denote the category of unbounded $\Aaa$-valued cochain complexes $A$, with differential 
$d:A^n\rightarrow A^{n+1}$. In \cite[arxiv version]{Hovey} a complete construction is given of the 
{\em injective model structure} on ${\rm cpx}(\Aaa )$. The cofibrations are degreewise monomorphisms and the weak equivalences are quasiisomorphisms; fibrations
are defined by the lifting property with respect to trivial cofibrations i.e. the monomorphic quasiisomorphisms. 

Hovey, refering to \cite{AFH97}, points out
that the fibrant objects are the complexes whose component objects $A^n$ are injective in $\Aaa$, and which are $K$-injective in the sense of Spaltenstein \cite{Spaltenstein}.
The injective model structure generalizes the notion of injective resolutions, and is well adapted to the construction of right derived functors, 
see \cite[arxiv version]{Hovey}, Proposition 2.14. 

Denote by $L{\cpx}(\Aaa )$ the Dwyer-Kan simplicial localization of the category ${\cpx}(\Aaa )$ inverting quasiisomorphisms
\cite{DwyerKan}. Existence of the injective
model structure shows that $L{\cpx}(\Aaa )$ is equivalent to a large simplicial category in the same universe in which ${\cpx}(\Aaa )$ is a large category. That is
to say the morphism classes in $L{\cpx}(\Aaa )$ may be replaced by simplicial sets, and heretofore we assume this is done. On the other hand, let 
$DG{\cpx}(\Aaa )$ denote the differential graded category as constructed by Kapranov in 
\cite{Kapranov0} whose objects are the fibrant objects of the injective model structure on ${\cpx}(\Aaa )$
and whose morphism complexes are just the complexes $Hom^{\cdot}(A^{\cdot} , B^{\cdot})$.

Proposition \ref{dgsimpl} gives a natural equivalence of simplicial categories 
$$
L{\cpx}(\Aaa ) \cong \widetilde{DP}(DG{\cpx}(\Aaa )).
$$

Suppose $\Lambda$ is a quasicoherent sheaf of $\Oo_X$-algebras on $X$. The above then applies
to the abelian category $\Aaa := \mdl _{\rm qc}(\Lambda )$ of sheaves of modules over $\Lambda$ which are quasicoherent as sheaves of $\Oo _X$-modules. 
We obtain the category $\cpx _{\rm qc}(\Lambda ):= \cpx (\mdl _{\rm qc}(\Lambda ))$ of unbounded complexes of sheaves of $\Lambda$-modules which are quasicoherent
as $\Oo _X$-modules. In the injective model structure, cofibrations are the termwise monomorphisms, so all objects are cofibrant. The fibrant objects are 
complexes of injective $\Lambda$-modules which are $K$-injective in the sense of Spaltenstein \cite{Spaltenstein}. 
The fibrations are dimensionwise split surjections with fibrant kernel \cite[arxiv version]{Hovey}, Prop. 2.12.

If a left-derived functor is called for, we need a different model structure. 
Hovey constructs the {\em projective model structure} on the category of complexes of modules over a ring
in \cite{HoveyBook}. This covers the case when $X=\Spec (A)$ is
affine and the ring $L=\Lambda (X)$ corresponds to $\Lambda$. A quasicoherent sheaf of $\Lambda$-modules is the
same thing as an $L$-module, so $\cpx _{\rm qc}(\Lambda ) = \cpx (L-\mdl )$ has the model structure
given by \cite{HoveyBook}. Recall from there: the fibrations are the termwise surjective morphisms. The cofibrations are the
termwise split injections with cofibrant cokernel. The cofibrant objects have a partial characterization:
cofibrant objects are termwise projective $L$-modules, and any bounded above complex of projective $L$-modules is 
cofibrant. 
Note that the terminology ``bounded below'' used in the statement of this criterion in \cite{HoveyBook} has to
be replaced by ``bounded above'' because we denote complexes cohomologically with differential of degree $+1$. 

The localization $L\cpx _{\rm qc}(\Lambda )$ is a simplicial category. 
Similarly we obtain DG-categories of fibrant and cofibrant objects $DG\cpx _{\rm qc}(\Lambda )$, corresponding to the
projective or injective structures. These can be distinguished by superscripts, and are linked by a diagram
of quasiequivalences of the form \eqref{twoway} described in \S \ref{dgcglob}. Either one is related to the simplicial localization by 
Proposition \ref{dgsimpl}.

Suppose $L$ is a ring and $\Aaa$   
an AB5 abelian category. Suppose we have a functor ${\bf V} : \mod (L)\rightarrow \Aaa$ that preserves direct sums and quotients. It extends to a functor 
${\bf V} : \cpx (L)\rightarrow \cpx ( \Aaa )$ which preserves colimits. This functor is completely determined
by the object ${\bf V}(L)\in \Aaa$ together with its right $L$-module structure, by the formula
$$
{\bf V}(E) = {\bf V}(L)\otimes _L E.
$$
Here the tensor product is defined in the obvious way if $E$ is a free $L$-module. For an arbitrary module $E$ the tensor product is
obtained by choosing a presentation of $E$ as the cokernel of a map between free modules; the tensor product with $E$ is the cokernel of the
corresponding map between tensor products with free modules. 

\begin{lemma}
\mylabel{genadjoint}
A functor ${\bf V}$ as above has as right adjoint
$$
{\bf V}^{\ast}: B\mapsto \Hom_{\Aaa }({\bf V}(L),B)
$$
on the category of complexes. This forms a Quillen pair, that is ${\bf V}$ is
a left Quillen functor from the projective model structure
on $\cpx (L)$ to the injective model structure on $\cpx (\Aaa )$.
\end{lemma}
\begin{proof}
The object $L$ considered as left $L$-module, has a right action
of $L$. The $\Hom$ from here transforms this back to a left action, so $\Hom_{\Aaa }({\bf V}(L),B)\in \cpx (L)$ 
whenever $B\in \cpx (\Aaa )$. One can verify the adjunction property using the expression for ${\bf V}$ as tensor product
with ${\bf V}(L)$. The recent reference \cite{BohmBrzWis} looks relevant here. 

The functor ${\bf V}$ is right exact by assumption, and it automatically 
takes cofibrations in the projective model structure of $\cpx (L)$ to cofibrations in the injective model structure of $\cpx (\Aaa )$ 
because cofibrations in the projective model structure are termwise split. We claim that ${\bf V}$ takes trivial cofibrations
to trivial cofibrations. If $E\rightarrow E'$ is a trivial cofibration then $E'/E$ is an acyclic cofibrant object.
We have an exact sequence in $\cpx (\Cc )$
$$
0\rightarrow {\bf V}(E)\rightarrow {\bf V}(E') \rightarrow {\bf V}(E'/E) \rightarrow 0.
$$
So it suffices to show that ${\bf V}(E'/E)$ is acyclic. Cofibrant objects of the projective model structure
are $K$-projective in the sense of Spaltenstein \cite{Spaltenstein}. This means that for any acyclic complex $A$,
the internal $Hom$ complex $\Hom (E'/E, A)$ is acyclic. This can be applied to $A=E'/E$ itself, and the identity
endomorphism is a closed element of degree zero. Hence it is the differential of an element of degree $-1$,
that is to say that the identity endomorphism of $E'/E$ is homotopic to zero. The same is therefore true of 
${\bf V}(E'/E)$, so ${\bf V}(E'/E)$ is acyclic. This shows that ${\bf V}(E)\rightarrow {\bf V}(E')$ is a trivial cofibration,
so ${\bf V}$ is a left Quillen functor. 
\end{proof}

\subsection{Complexes of sheaves on a site}
\mylabel{sec-cpx}

We collect here some observations and definitions concerning the model categories of
unbounded complexes of sheaves of $\Oo$-modules on a site $\Xx$. If $Y\in \Xx$ the category $\Xx /Y$ of objects over $Y$
has a natural induced Grothendieck topology.

Suppose $U$ and $V$ are complexes of sheaves of $\Oo$-modules on a site $\Xx$. Then
we can define the {\em internal $Hom$ complex} $\uHom^{\cdot} (U,V)$ which is the complex of sheaves 
$\uHom^{\cdot} (U,V)(x):= \Hom^{\cdot}(U|_{\Xx /x}, V|_{\Xx /x} )$.

Suppose we are given a presheaf of sets $Y$ on $\Xx$, and a sheaf of $\Oo$-modules $U$ on $\Xx /Y$.
Let $p:Y\rightarrow \ast$ denote the projection to the trivial punctual sheaf. For any object $W\in\Xx$, we can form the $\Oo (W)$-module 
$\Hom_{\Xx /Y}(W|_{\Xx /Y}, U)$. As $W$ varies this is a sheaf on $\Xx$ which we denote by $p_{\ast}(U)$.

\begin{lemma}
\mylabel{adjunction}
The construction $p_{\ast}$ is right adjoint to the restriction, and that extends to complexes of
sheaves by applying it termwise to the sheaves in each complex:
$$
\uHom ^{\cdot}(U|_{\Xx / Y}, V) = \uHom ^{\cdot} (U,p_{\ast}(V)).
$$
\end{lemma}
\eop

We say that a complex of sheaves of $\Oo$-modules $A$ on a ringed site $\Xx$  is {\em cohomologically flasque} if, for any fibrant resolution
$A\rightarrow A'$, and for any object $Y\in \Xx$, the map $A(Y)\rightarrow A'(Y)$ is a quasiisomorphism of 
complexes of abelian groups. This is similar to Jardine's notion of flasque simplicial presheaf \cite{Jardine}, and is an analogue
of the condition of being an $n$-stack or Segal $n$-stack considered in \cite{HirschowitzSimpson}.

\begin{lemma}
\mylabel{cohflasque1}
A fibrant complex $A$ is cohomologically flasque. Any quasiisomorphism between cohomologically flasque complexes
$U\rightarrow V$  induces a quasiisomorphism on each $U(X)\rightarrow V(X)$. We can test cohomological flasqueness of a complex $B$ by looking
at any single fibrant resolution $B\rightarrow B'$ instead of all possible ones. If $\Aff$ is the site of affine schemes of finite type over $\cc$
with the etale topology, a bounded complex on $\Aff /X$ whose elements are locally free $\Oo$-modules
(possibly of infinite rank) is cohomologically flasque.
\end{lemma}
\eop

\begin{theorem}
\mylabel{cohflares}
Suppose $Y$ is a presheaf of sets over $\Xx$. If $A$ is a fibrant (resp. cohomologically flasque) complex
of sheaves of $\Oo$-modules on $\Xx$, then $A|_{\Xx /Y}$ is a fibrant (resp. cohomologically flasque) complex
of sheaves of $\Oo$-modules on $\Xx /Y$.
\end{theorem}
\begin{proof} 
Consider the functor $p:\Xx /Y \rightarrow \Xx$. The restriction $p^{\ast}$ has as left adjoint the functor
defined as follows: first define  $p_{!}^{\rm pre}$ from presheaves of $\Oo$-modules on $\Xx /Y$
to presheaves of $\Oo$-modules by 
$$
p_{!}^{\rm pre}(A)(Z):= \bigoplus _{f:Z\rightarrow Y}A(Z, f);
$$
Then let $p_{!}(A)$ denote the sheaf of $\Oo$-modules associated to $p_{!}^{\rm pre}(A)$. 
These form a Quillen pair, so the restriction $p^{\ast}$ of a fibrant object is fibrant; hence the same for cohomological flasqueness. 
\end{proof}

The assignment to $x\in \Xx$ of the injective model category of unbounded complexes $\cpx ^{\rm inj}(\Xx /x, \Oo )$ is a left Quillen presheaf. We get a
prestack $[x\mapsto L\cpx (\Xx /x, \Oo )]$ over $\Xx$. 

\begin{lemma}
\mylabel{dpfibrant}
Suppose $A$ and $B$ are objects in the injective model category ${\rm cpx}(\Xx , \Oo )$, with $B$ fibrant. Then 
$[ x\mapsto DP\tau_{\leq 0}\uHom ^{\cdot}(A,B)(x)]$ is
a fibrant simplicial presheaf on $\Xx$. Hence, 
the prestack $[x\mapsto L\cpx (\Xx /x, \Oo )]$ is a protostack in the notation of \cite[\S 15]{HirschowitzSimpson}.
\end{lemma}
\begin{proof}
Suppose $g :U\rightarrow V$ is a trivial cofibration of simplicial presheaves on $\Xx$. 
The lifting property for $[ x\mapsto DP\tau_{\leq 0}\uHom ^{\cdot}(A,B)(x)]$ along $g$
corresponds to the lifting property for $\uHom^{\cdot}(A,B)$ along $N(\cc g) : N(\cc U) \rightarrow N(\cc V)$,
and $N(\cc g)$ is again a trivial cofibration. 
\end{proof}

Recall from \cite[\S 15]{HirschowitzSimpson} 
that the condition ``sufficiently many disjoint sums'' means that there is an ordinal $\beta$ such that the covering families of
cardinality $\gamma < \beta$ generate the the topology, and that disjoint sums of any cardinality $\gamma <\beta$ exist.
The descent result \cite[Th\'eor\`eme 19.4]{HirschowitzSimpson} applied to 
the left Quillen presheaf of model categories  $x\mapsto {\rm cpx}(\Xx , \Oo )$ can be strengthened. 
We sketch the proof which was obtained in conversations with A. Hirschowitz.

\begin{theorem}
\mylabel{like194}
Suppose $(\Xx , \Oo )$ is a ringed site, and assume one of the following conditions:
\newline
(1)\,  $\Xx$ admits fiber products and sufficiently many disjoint sums; or
\newline
(2)\, the representable presheaves on $\Xx$ are sheaves, and our universe contains an inaccessible cardinal. 
\newline
Then $[x\mapsto L\cpx (\Xx /x, \Oo )]$
is a stack over $\Xx$. Furthermore, the $(\infty ,1)$-category of global sections---defined by first taking a 
fibrant replacement---is naturally equivalent to the $(\infty ,1)$ category of global complexes $L\cpx  (\Xx , \Oo )$. 
\end{theorem}
\begin{proof}
(1) The arguments of 
the proof of \cite[Corollaire 21.1]{HirschowitzSimpson} to check conditions 
conditions (0)-(4)  of \cite[Th\'eor\`eme 19.4]{HirschowitzSimpson} work equally well for unbounded complexes, once we know about the injective model structure, and Lemma \ref{dpfibrant} gives hypothesis (5). 

\noindent
(2)\,  Let $\uu$ be an intermediate universe. Apply (1) to
the site ${\bf Sh}(\Xx , \uu )$  of $\uu$-sheaves of sets on $\Xx$ with the canonical topology. Representables give a
functor $h:\Xx \rightarrow {\bf Sh}(\Xx , \uu )$. 
For each object $x$ the categories
of sheaves on $\Xx /x$ and ${\bf Sh}(\Xx , \uu )/x$ are equivalent, so 
$L\cpx (\Xx /x,\Oo ) \cong L\cpx ({\bf Sh}(\Xx , \uu )/x, \Oo )$. 
One can show that 
the functor $h_{!}$ preserves $\uu$-small hypercovers; therefore $h^{\ast}$ preserves the $(\infty ,0)$-stack condition \cite{DuggerHollanderIsaksen}, and the argument 
of \cite[Lemme 10.5]{HirschowitzSimpson}
shows that $h^{\ast}$ preserves the $(\infty ,n)$-stack condition. Thus $[x\mapsto L\cpx (\Xx /x,\Oo )]$ is an $(\infty ,1)$-stack over $\Xx$.

For the last part of the theorem, do the whole theory using sheaf objects 
as in Joyal's original
\cite{Joyal}, rather than presheaf objects. In this context $h_{!}$ preserves cofibrations and trivial cofibrations (which is not the case if we use presheaves).
With this, the restriction of the fibrant replacement over ${\bf Sh}(\Xx , \uu )$ is a fibrant replacement over $\Xx$ and we get the
required identification. 
\end{proof}

\subsection{Perfect complexes}
\mylabel{sec-perfect}

Assume now that $\Xx = {\Aff}$ is the site of affine schemes of finite type over $\cc$ with the etale topology and ring
$\Oo$. Suppose $X=\Spec (A) \in \Aff$. If $U$ is an $A$-module corresponding to a
quasicoherent sheaf $\tilde{U}$ on $X_{\rm zar}$, then we obtain a ``quasicoherent'' sheaf of $\Oo$-modules $\tilde{U}(\Spec  (B)):= U\otimes _AB$ on $\Aff /X$.
On a site of the form $\Aff /\Gg$, the {\em quasicoherent} sheaves of $\Oo$-modules are those which are quasicoherent in the above sense 
when restricted to any $\Aff /X$.

Suppose $a\leq b$ are integers.
For $X\in \Aff$, a complex of sheaves of $\Oo$-modules $E$ on $\Aff /X$ is {\em strictly perfect of amplitude $[a,b]$} if $E$ is isomorphic to a complex
of quasicoherent sheaves whose components $E^i$ are locally free of finite rank and nonzero only in the interval $a\leq i \leq b$. 
A complex of sheaves $E$ of $\Oo$-modules on $\Aff /X$ is {\em perfect of amplitude $[a,b]$} if there is a covering of $X$ by Zariski open sets
$Y_j$ such that the restriction of $E$ to each $\Aff /Y_j$ is quasiisomorphic, as a complex of sheaves of $\Oo$-modules, to a strictly perfect
complex of amplitude $[a,b]$. We say that $E$ is {\em perfect} if it is Zariski-locally perfect of some amplitude. 
The amplitude can be chosen globally on $X$ by quasicompactness. 

Let $\Perf$ denote the stack which to $X\in \Aff$ associates the full $(\infty , 1)$-subcategory 
$$
\Perf (X)\subset L\cpx (\Aff /X, \Oo )
$$
consisting of perfect complexes. For a fixed interval of amplitude $[a,b]$ let
$\Perf ^{[a,b]}\subset \Perf$ be the full substack of complexes which are perfect of amplitude $[a,b]$. Descent for the property of being perfect
means that $\Perf$ and $\Perf ^{[a,b]}$ are substacks of $X\mapsto \cpx (\Aff /X, \Oo )$.

A. Hirschowitz suggested applying the notion of $n$-stack to the moduli of perfect complexes. This suggestion, which 
led to \cite{HirschowitzSimpson}, was motivated by his study
\cite{Hirschowitz} where the utility of working on a big site became clear. 

We recall here the result of Toen and Vaqui\'e on the local geometricity of the moduli stack of perfect complexes \cite{ToenVaquie}.
Their theorem is the main ingredient in our argument.

\begin{theorem}[Toen-Vaqui\'e \cite{ToenVaquie}]
\mylabel{tv}
Suppose $X\rightarrow S$ is a smooth projective morphism of relative dimension $d$ polarized by $\Oo _X(1)$. 
The stack $\Perf (X/S)$ is locally geometric. It is covered by geometric $n$-stacks of the form 
$\Perf ^{[a,b]}(\underline{h})$ where $\underline{h}(i,n)$ is a function defined for
$a\leq i\leq b+d$ and  $n$ in some interval $M\leq n \leq N$ depending on $X/S$. For $Y\rightarrow S$ the objects $E\in \Perf ^{[a,b]}(\underline{h})(Y)$
are the perfect complexes on $X\times _SY$ such that
${\rm dim} \hh ^i(X_y, E (n)|_{X_y})\leq \underline{h}(i,n)$ for $a\leq i \leq b+d$ and 
$M\leq n \leq N$.
\end{theorem}

If $\Lambda$ is a quasicoherent sheaf of $\Oo _X$-algebras, then 
in either the projective or injective model structures, consider the full subcategories denoted 
$L\cpx _ {\Oo -{\rm perf}}(\Lambda )$ or $DG\cpx _ {\Oo -{\rm perf}}(\Lambda )$ of complexes of sheaves of $\Lambda$-modules whose components are $\Oo _X$-quasicoherent, and which are perfect
complexes of $\Oo _X$-modules, that is to say locally quasiequivalent to bounded complexes of free $\Oo _X$ modules of finite rank. 
The Dold-Puppe compatibility says that
$$
L\cpx _{\rm qc}(\Lambda ) = \widetilde{DP}(DG\cpx _{\rm qc}(\Lambda )), \;\;\; 
L\cpx _{\Oo -{\rm perf}}(\Lambda ) = \widetilde{DP}(DG\cpx _{\Oo -{\rm perf}}(\Lambda )).
$$

In the next section we will look at a $1$-stack of groupoids $\Gg$ over $\Aff$, provided also with
a morphism from a scheme $X\rightarrow \Gg$. 
A complex of sheaves of $\Oo$-modules on $\Aff /\Gg$, which we should no longer require to be quasicoherent, is {\em $\Oo$-perfect}
if its restriction to $\Aff /X$ is quasiisomorphic to a perfect complex on $X$. Again we can define full subcategories of 
the $(\infty ,1)$ or dg categories of complexes denoted by a similar notation.


\section{Complexes on the crystalline site associated to a formal category}
\mylabel{strictif-crys}

\subsection{Formal categories}
\mylabel{fc}

A formal category of smooth type $(X,\Ff )$ consists of a separated scheme of finite type $X$ and 
a formal scheme $\Ff$ with structure maps
$$
X\stackrel{i}{\rightarrow} \Ff \tworightarrows X, \;\;\;
\Ff \times _ X \Ff \stackrel{m}{\rightarrow} \Ff
$$
inducing a functor from schemes to  groupoids, such that $\Ff$ is supported along the closed subscheme $i(X)$
and satisfies a local freeness property recalled below.

Denote by $\Gg := [X/\! / \Ff ]$ the $1$-stack on $\Aff$ 
associated to the prestack of groupoids defined by $(X,\Ff )$, which may also be considered
as any kind of $\infty$-stack. The formal category can be recovered from the map $X\rightarrow \Gg$. 

Corresponding to $\Ff$ is a complete sheaf of algebras ${\bf f}$ over $\Oo _{X\times X}$ supported along
the diagonal. 
Let ${\bf j}\subset {\bf f}$ denote the ideal of
$i(X)$. Considered as a sheaf on $X$, ${\bf f}$ has
a structure of both left and right $\Oo _X$-algebra, and it is ${\bf j}$-adically complete with ${\bf j}$ being a two-sided ideal. 
The {\em smooth type} condition means that ${\bf j}^k/{\bf j}^{k-1}$ are locally free $\Oo _X$-modules of finite rank and
$$
Gr_{{\bf j}}({\bf f}) \cong Sym^{\cdot}(\Omega ^1_{X,\Ff }), \;\;\; \mbox{where}\;\;\;
\Omega ^1_{X,\Ff }:= {\bf j}/ {\bf j}^2.
$$
Define the corresponding ``tangent sheaf'' by
$$ 
\Theta _{X,\Ff}:= \underline{\Hom} _{\Oo _X}(\Omega ^1_{X,\Ff},\Oo _X). 
$$

We can look at the {\em big crystalline site} $CRIS(X,\Ff ) := \Aff / \Gg $ whose objects are pairs
$(Y,g)$ with $Y\in \Aff$ and $g\in \Gg (Y)$. 
The site is ringed by the pullback of the structure sheaf from $\Aff$, denoted here also by 
$\Oo$. The morphism $X\rightarrow \Gg$ corresponds to a sheaf of sets still denoted by $X$, and
$CRIS(X,\Ff )/X \cong \Aff /X$. 

A complex of sheaves of $\Oo$-modules on $CRIS(X, \Ff )$ is called $\Oo$-{\em perfect} if its restriction to 
$\Aff /X$ is perfect in the sense of \S \ref{sec-perfect}.

If $X$ is a smooth scheme and $\Ff _{DR} = (X\times X)^{\wedge}$ is the de Rham formal category, then
$CRIS(X,\Ff _{DR})$ is equivalent to a big crystalline site of $X$ whose objects are pairs $(U,Y)$
where $U\rightarrow X$ is a reduced affine scheme of finite type mapping to $X$ and $U\subset Y$ is an infinitesimal thickening,
with the etale topology. This is basically the same as the usual big crystalline site considered in \cite{Berthelot} 
\cite{BerthelotOgus} \cite{BerthelotIllusie}; the extra degree of freedom in the usual case where $U$ can be non-reduced,
doesn't seem to be essential. 

The following theorem relates complexes on $CRIS(X, \Ff )$ to 
$Hom$ in the $(\infty ,2)$-category of $(\infty ,1)$-stacks over $\Aff$, from $\Gg$ to the
stack of complexes. It basically reflects the principle that stacks over
$CRIS(X, \Ff )$ are the same  thing as stacks over $\Aff$ together with maps to $\Gg$. 
Use the Segal approach of \cite{HirschowitzSimpson} to make these notions precise. 

\begin{theorem}
\mylabel{cpxonstack}
If $(X,\Ff )$ is a formal category of smooth type,
there is a natural equivalence of $(\infty , 1)$-categories
$$
\Hom _{(\infty , 1){\rm St}/\Aff }(\Gg , [x\mapsto L\cpx (\Aff /x,\Oo )])
\stackrel{\sim}{\rightarrow}
L\cpx (CRIS(X, \Ff ),\Oo ).
$$
The full subcategory of $L\cpx (CRIS(X, \Ff ),\Oo )$ consisting of perfect complexes is naturally equivalent to 
$\Hom _{(\infty , 1){\rm St}/\Aff }(\Gg , \Perf )$. 
\end{theorem}
\begin{proof}
Let $\Cc$ denote the stack $[x\mapsto L\cpx (\Aff /x,\Oo )]$ on $\Aff$. Note that
the over-categories are the same $CRIS(X,\Ff )/(Y,g) \cong \Aff /Y$ and they share the same hypercovers, 
so the pullback of $\Cc$ to  $CRIS(X,\Ff )$ is again
a stack. 

Let $\Cc '$ be a fibrant replacement for $\Cc$ in the injective model category \cite{HirschowitzSimpson}
of Segal stacks over $\Aff$. 
Let $(\Cc '|_{CRIS(X,\Ff )})'$ be a fibrant replacement over $CRIS(X,\Ff )$. 
There is  $\tau \in \Gamma ( CRIS(X,\Ff ) ,  \Gg |_{CRIS(X,\Ff )})$ a tautological section
which to $(Y,g)$ associates $g$. Composition with $\tau$ induces the morphism in question 
using
the last sentence of Theorem \ref{like194}
$$
\Hom (\Gg , \Cc ')\rightarrow \Gamma ( CRIS(X,\Ff ) ,(\Cc '|_{CRIS(X,\Ff )})')\approx L\cpx (CRIS(X,\Ff ),\Oo ).
$$

To study sections of the fibrant replacement of the pullback, we use a formulation in terms of fibered categories.
This theory is developped in the world of quasicategories by Lurie in \cite[\S 2.3]{LurieTopos} and can be adapted to the
present context. 
In particular, corresponding to the
stack $\Cc$ over $\Aff$ there is a {\em fibered $(\infty , 1)$-category} denoted $\int _{\Aff} \Cc \rightarrow \Aff$. 
Note that $\int _{\Aff}\Gg = \Aff /\Gg = CRIS(X,\Ff )$, these are just different notations for the same thing. The section $\tau$ is the diagonal
in this interpretation. As in the usual $1$-categorical case, there is a notion of {\em cartesian morphism}
in $\int _{\Aff}\Cc$, which is an arrow $(x,f)\rightarrow (y,g)$ covering an arrow $a:x\rightarrow y$, such that
the corresponding map $x\rightarrow a^{\ast}(g)$ is an equivalence in the $(\infty , 1)$-category $\Ff (x)$.  

Before trying to take the $(\infty ,1)$-category of cartesian sections or cartesian morphisms, one should 
make a fibrant replacement of the fibered category in the model category of $(\infty ,1)$-categories;
call this the ``fibrant integral'' $\left( \int _{\Aff}\Cc \right) '\rightarrow \Aff$.   
Then $\Gamma _{(\infty , 1)}(\Aff , \Cc ' )$ is equivalent to the $(\infty , 1)$-category of cartesian sections of the fibrant integral
of $\Cc$ i.e. those which send pullback maps to equivalences. The same holds true for $(\Cc ' |_{\Aff /\Gg})'$
over $\Aff /\Gg$. However, the advantage of this point of view is that the pullback from $\Aff$ to $\Aff /\Gg$ of the fibrant integral 
of $\Cc$, is already fibrant over $\Aff /\Gg$, in other words
$\left( \int _{\Aff /\Gg}\Cc |_{\Aff / \Gg  } \right) '=  \left( \int _{\Aff}\Cc \right) '\times _{\Aff}\Aff / \Gg $. 
A section  from $\Aff / \Gg $ into here is cartesian 
if and only if its first projection is cartesian 
as a morphism $\Aff / \Gg = \int _{\Aff}\Gg \rightarrow\left( \int _{\Aff}\Cc \right) '$ of fibered categories over $\Aff$,
because all morphisms of $\Aff / \Gg$ are cartesian: $\Gg$ is a stack of groupoids.

On the other hand, the $(\infty , 1)$-category of these cartesian morphisms is equivalent to
$\Hom _{(\infty , 1){\rm St}/\Xx } (\Gg , \Cc ' )$.
After making the above reinterpretations, the map in question is this equivalence. 
\end{proof}

\subsection{Morphisms}
\mylabel{morph}

Suppose now that $(X, \Ff )$ and $(Y, \Hh )$ are formal categories of smooth type. A {\em morphism of formal categories} $f$ is a pair of morphisms 
$f_{\rm ob} : X\rightarrow Y$ and $f_{\rm mor}:\Ff \rightarrow \Hh$ compatible with the categorical structure maps. 
Such a morphism induces a morphism of sheaves of $\Oo _X$-modules on $X$, 
$$
df: \Theta _{X,\Ff} \rightarrow f_{\rm ob}^{\ast} \Theta _{Y,\Hh}.
$$
We say that $f$ is {\em a morphism
of smooth type} if both the source and target are formal categories of smooth type,
if $f_{\rm ob}:X\rightarrow Y$ is flat, and if $df$ is surjective.

A morphism $f$ gives $[X/\! /\Ff ]\rightarrow [Y /\! /\Hh ]$ hence a functor of crystalline sites 
$$
\{ f \} : CRIS(X,\Ff )\rightarrow CRIS (Y, \Gg ). 
$$
This morphism gives a higher direct image functor
of $(\infty , 1)$-categories
$$
R\{ f \} _{\ast} : L\cpx (CRIS(X, \Ff ), \Oo )\rightarrow L\cpx (CRIS(Y, \Gg ), \Oo ).
$$
On the other hand, pullback is a functor
$$
\{ f \} ^{\ast} :  L\cpx (CRIS(Y,\Gg  ), \Oo )\rightarrow  L\cpx (CRIS(X,\Ff ), \Oo ).
$$

\begin{lemma}
\mylabel{fcdirim}
If $U$ is a perfect complex on $CRIS(Y, \Gg )$ then the pullback $\{ f \} ^{\ast}(U)$ is a perfect complex on $CRIS(X,\Ff )$.
If $f$ is a morphism of smooth type and 
$U$ is a perfect complex on $CRIS(X,\Ff )$ then  $\rr \{ f \} _{\ast} (U)$ is a perfect complex on $CRIS(Y, \Gg)$.
\end{lemma}
{\em Proof:}
For pullback, this follows from the fact that the pullback of a perfect complex from $\Aff /Y$ to $\Aff /X$ remains perfect.
For higher direct image, the argument is a standard one---see for example Katz \cite{Katz}---which we don't reproduce here. Many elements going into this
argument, such as the construction of a de Rham-Koszul resolution, will occur in the subsequent discussion below. 
\eop

\subsection{The sheaf of rings of differential operators}
\mylabel{diffop}

Classically \cite{Berthelot}
associated to a formal category of smooth type $(X,\Ff )$  we get a filtered sheaf of rings of differential operators $\Lambda$ defined by
$$
\Lambda  = \bigcup J_k\Lambda,
\;\;\; J_k \Lambda := \Hom _{\Oo _X, {\rm left}}({\bf f}/{\bf j}^k, \Oo _X).
$$
We can think of $\Lambda$ as the sheaf of morphisms from ${\bf f}$ to $\Oo _X$ which locally factor through some ${\bf f}/{\bf j}^k$. 
The composition operation of $\Ff$ gives the ring structure making $\Lambda$
into a split almost-polynomial sheaf of rings of differential  operators in the sense of \cite{moduli1}. 
Dual to
the first projection $\Ff \rightarrow X$ we get
$J_1\Lambda \rightarrow \Oo _X$ which gives the splitting $J_1\Lambda \cong \Theta _{(X,\Ff )}\oplus \Oo _X$. 

In much of the discussion below we will localize to assume that $X=\Spec  (A)$ is affine. The formal category $\Ff$ 
corresponds to a complete
$A\otimes _{\cc}A$ algebra $F:= {\bf f}(X\times X)$ considered as a left and right $A$ 
algebra. Then $\Lambda$ is determined by its ring of global sections
$L:= \Lambda (X)$ with filtration by $J_kL$ and 
identification $J_0L =A$,
in particular $L$ also has structures of left and right $A$-algebra. 

Introduce the following notation: the sign $<$ 
corresponds to the
left structure, and the sign $>$ corresponds to the right structure	. These 
are included in the
subscript, so for example $U\otimes _{<A<} V$ would mean the tensor product 
of $U$ and $V$ using the left structures
on both sides. The default is $U\otimes _{>A<}V$ if that makes sense, otherwise 
$U\otimes _{<A<}V$. 
When necessary adopt a similar convention for $Hom$.

Dual to the ring multiplication operation of $F$ is a coalgebra structure on 
$L$. Note that
$$
L\otimes _{<A<}L = Hom _{<,<,A}(F\otimes _{<  A <}F, A).
$$
The ring multiplication $r:F\otimes _{<  A <}F\rightarrow F$ thus gives a map 
of left $A$-modules
$$
\Delta : L\rightarrow L\otimes _{<A<}L.
$$
The expression $\Delta (\varphi ) = \sum _j \psi _j\otimes \eta _j$ is determined as the unique one such that
$\varphi (r(f\otimes g)) = \sum _j \psi _j(f)\eta _j(g)$.
Calculation of the coproduct in the case $\Lambda = \Dd _X$ leads to the introduction of ``shuffles''. 

The map $\Delta$ actually goes into a special submodule of $L\otimes _{<A<}
L$. This tensor product, contracting the two left-module structures,
has two distinct structures of right $A$-module,
coming respectively from the first and second factors. Let 
$(L\otimes _{<A<}L)^{=}$ denote the subgroup of tensors $\sum _j\psi 
_j\otimes \eta _j$ with the
property that for any $a\in A$, there is equality 
$$
\sum _j(\psi _ja)\otimes \eta _j = \sum _j\psi _j\otimes (\eta _j a). 
$$
Then $(L\otimes _{<A<}L)^{=}$ has a unique right $A$-module structure.

The coproduct $\Delta$ takes image in $(L\otimes _{<A<}L)^{=}$ 
and intertwines the right $A$-module structure of $L$ with
the right module structure on $(L\otimes _{<A<}L)^{=}$. 

We can now define the tensor product of $L$-modules which will later correspond to 
tensoring sheaves on $CRIS (X,\Ff )$. Suppose $U$ and $V$
are two left $L$-modules. The module structures correspond to maps of left $L$-modules
$$
L\otimes _{>A<}U\rightarrow U, \;\;\; L\otimes _{>A<}V\rightarrow V.
$$
The left module structure on the domain comes from that of $L$. These two 
structures combine to give a map of left $A$-modules 
$$
(L\otimes _{<A<}L)^{=}\otimes _{>A<}(U\otimes _{<A<}V)\rightarrow 
(U\otimes _{<A<}V).
$$
If $\sum _j\psi _j\otimes \eta _j\in (L\otimes _{<A<}L)^{=}$ and $u\in U$ and 
$v\in V$
then
$$
(\sum _j\psi _j\otimes \eta _j)\otimes (u\otimes v)\mapsto \sum _j \psi _j(u)
\otimes \eta _j(v).
$$
This is compatible with the rule $u\otimes (av) = (au)\otimes v$ because of the 
condition 
$\sum _j\psi _j\otimes \eta _j\in (L\otimes _{<A<}L)^{=}$. It is also clearly 
a map of left $A$-modules. 

Now compose this map with $\Delta : L\rightarrow (L\otimes _{<A<}L)^{=}$ 
(which is a map of right $A$-modules)
to get a map
$$
L\otimes _{>A<}(U\otimes _{<A<}V)\rightarrow U\otimes _{<A<}V.
$$
This is the $L$-module structure on the tensor product. Note that it has to 
involve the shuffling operation inherent in $\Delta$. 
We leave to the reader to verify that it is indeed a structure of left $L$-
module, using the coassociativity of $\Delta$. Denote the tensor product 
with this $L$-module structure by $U\otimes ^{\Delta}_{<A<}V$.

\begin{lemma}
\mylabel{interchange}
Suppose $U$ is a left $L$-module. Consider $L$ itself as a left $L$-module. 
Then there is a natural isomorphism
$$
L\otimes _{<A<}U \cong L\otimes _{>A<}U
$$
taking 
the left $L$-module structure $L\otimes ^{\Delta}_{<A<}U$ defined above, to the tautological left module structure
on $L\otimes _{>A<}U$ which uses only the structure on $L$ and not on $U$. 
\end{lemma}
\begin{proof}
Let $>1$ and $>2$ denote the two different right $A$-module structures of 
$L\otimes _{<A<}L$, both preserved by $\Delta$. Tensoring with $U$ we 
get a map
$$
L\otimes _{>A<}U \rightarrow (L\otimes _{<A<}L)\otimes _{>2,A<}U.
$$
Composed with the left $L$-module structure $L\otimes _{>A<} U \rightarrow U$
we get the map 
$$
L\otimes _{>A<}U \rightarrow (L\otimes _{<A<}L)\otimes _{>2,A<}U\cong 
L\otimes _{<A<}(L\otimes _{>A<} U)
\rightarrow L\otimes _{<A<}U
$$
which one can check is an isomorphism. 
\end{proof}

\subsection{The de Rham-Koszul resolution}
\mylabel{koszul}

In the classical theory of $\Dd_X$-modules, we have a resolution of the structure sheaf $\Oo$ considered as a $\Dd _X$-module
corresponding to the de Rham complex. This extends to the case of a formal category of smooth type, as I learned from Berthelot \cite{Berthelot}.
His assumption that $\Omega ^1_{(X,\Ff )}$ be generated by differentials is unnecessary, see also \cite[Ch. VIII]{Illusie}.
The existence of a $\Lambda$-module structure on $\Oo$ is a consequence of the fact that 
$\Lambda$ comes from a formal category of smooth type: there can be almost-polynomial sheaves of rings of differential operators
not having such a module. Given that $\Lambda$ corresponds to $(X, \Ff )$, the structure sheaf $\Oo$ on the site $\Aff /\Ff$ corresponds
to a $\Lambda$-module whose underlying quasicoherent sheaf on $X$ is $\Oo _X$. 

We have a de Rham-Koszul resolution
of left $\Lambda$-modules 
$$
\Kk ^{\cdot}:= \left(  \Lambda \otimes _{\Oo} \bigwedge ^r \Theta _{(X,\Ff )}\rightarrow \ldots \rightarrow \Lambda \otimes _{\Oo}\Theta  _{(X,\Ff )}
\rightarrow \Lambda \right) \rightarrow \Oo . 
$$
If $E$ is a complex of left $\Lambda$-modules then 
tensoring with $\Kk ^{\cdot}$ using the tensor product which we defined above using the coproduct $\Delta$
yields a resolution by locally free left $\Lambda$-modules 
$$
\Kk ^{\cdot}\otimes ^{\Delta}_{<\Oo _X<} E \rightarrow E. 
$$

\begin{proposition}
\mylabel{koszulcor}
An $\Oo _X$-perfect complex of quasicoherent left $\Lambda$-modules $E$ is quasiisomorphic, locally over $X$,
to a bounded complex of finite rank locally free left $\Lambda$-modules. 
\end{proposition}
\begin{proof}
Use the tensored de Rham-Koszul resolution above. 
Restrict to an open covering on which the terms of $\Kk ^{\cdot}$ are free $\Lambda$-modules. Here the interchange property Lemma \ref{interchange}
for the tensor product allows us to exchange the tensor product $\otimes ^{\Delta}_{<\Oo _X<}$ with the easier $\otimes _{>\Oo _X<}$. 
If $E$ is $\Oo _X$-perfect, choose a quasiisomorphism of complexes of $\Oo _X$-modules between $E$ and a strictly perfect complex $E'$. 
The tensor product $\otimes _{>\Oo _X<}$ transforms this quasiisomorphism into a quasiisomorphism of complexes of left $\Lambda$-modules,
and since $E'$ is strictly perfect, $\Kk ^{\cdot}\otimes _{>\Oo _X<} E '$ is a bounded complex of finite rank free $\Lambda$-modules. 
\end{proof}

\subsection{Crystalization}
\mylabel{sec-crys}
In this section we consider in some detail the functor going from $L$-modules to sheaves on $CRIS(X,\Ff )$. 
Start with the crystalization of modules.
Following Berthelot \cite{Berthelot} we can interpret $\Oo$-quasicoherent sheaves on $CRIS(X,\Ff )$ as $\Oo _X$-quasicoherent $\Lambda$-modules.
In  the affine case $X=\Spec  (A)$ $L=\Lambda (X)$, an $\Oo$-quasicoherent sheaf on 
$CRIS(X,\Ff )$ is given by an $A$-module $V$, together with
an isomorphism 
$$
\xi : F\otimes _{>A<}U \cong F\otimes _{<A<}U.
$$
It is required to satisfy a cocycle condition \cite{Berthelot}.
The cocycle condition can be generalized to the case of stacks or higher stacks, which is one of the
subjects of the papers of Breen and Messing \cite{Breen} \cite{BreenMessing1} \cite{BreenMessing2}. 
We get around the manipulation of such higher cocycles by using a model-category approach, but of course 
the cocycles represent what is really going on. 

The isomorphism $\xi$ corresponds to an action of $L$ which is dual to $F$. 
In the direction going from $\Lambda$-modules  to sheaves on $CRIS(X,\Ff )$, we denote this construction by $\Diamond$. 
If $U$ is a left $L$-module then $\Diamond (U)$ on $CRIS(X,\Ff )$ corresponds to the same underlying
$A$-module $U$, with a multiplication structure
$\xi_1 : U\rightarrow F\otimes _{<A<}U$
coming from the $L$-module structure of $U$. This extends by $F$-linearity to a map $\xi$ satisfying the cocycle condition because of associativity of the action
of $L$, see \cite{Berthelot}.  
The pair $(U,\xi )$ corresponds to a sheaf of $\Oo$-modules denoted $\Diamond (U)$ on $CRIS(X,\Ff )$.

The tensor product defined in \S \ref{diffop} using $\Delta$ corresponds via $\Diamond$ to the usual tensor product 
on $CRIS(X,\Ff )$, that is $\Diamond (U\otimes ^{\Delta}_{<A<}V) = \Diamond (U)\otimes _{\Oo }\Diamond (V)$.

The construction $\Diamond$ extends to complexes of $L$-modules, and can even be defined in the global non-affine situation
from $\Oo_X$-quasicoherent complexes of $\Lambda$-modules to complexes on  $CRIS(X,\Ff )$. However, in the global case we cannot derive $\Diamond$
directly because there is no projective model category structure on complexes of $\Lambda$-modules. The globalization of $\Diamond$ will
be treated in \S \ref{sec-lambdamodel} below using the notion of $(\infty , 1)$-stack. 

Getting back to the affine case, we are in exactly the situation of Lemma \ref{genadjoint} with
 $\Aaa = \mdl (CRIS(X, \Ff ),\Oo  )$ and the functor ${\bf V} = \Diamond$.
We obtain a left Quillen functor from the projective model category of $L$-modules to the injective model category
of sheaves of $\Oo$-modules on $CRIS(X, \Ff )$
$$
\Diamond : \cpx ^{\rm proj}(L) \rightarrow \cpx ^{\rm inj}(CRIS (X,\Ff ),\Oo )
$$
which derives to give the functor between $(\infty , 1)$-categories
$$
\Diamond : L\cpx (L )\rightarrow L\cpx (CRIS(X,\Ff ),\Oo ).
$$
The right Quillen adjoint $\Diamond ^{\ast}$ will be used in the next section to strictify complexes on the crystalline site.

\subsection{Strictification of sheaves on the formal category}
\mylabel{sec-from}

If $(X,\Ff )$ is a formal category with $X=\Spec (A)$ affine, we show that the crystallization functor
$\Diamond$ defined in the previous section is an equivalence on categories of $\Oo$-perfect objects. 
This argument is inspired by \cite[\S 5]{BerthelotOgus}. 

Consider the functor $p:\Aff /X\rightarrow CRIS(X, \Ff )$. 
If $B$ is a complex of $\Oo$-modules on $CRIS(X, \Ff )$ denote its restriction to $\Aff /X$ by $B|_X:= p^{\ast}(B)$. If $B$ is fibrant
(resp. cohomologically flasque \S \ref{sec-cpx}) then so is $B|_X$ by 
Theorem \ref{cohflares}.
On the other hand, if $A$ is a complex of $\Oo$-modules on ${\Aff}/X$
then we can form the complex $p_{\ast}A$ on $CRIS(X, \Ff )$.

Say that a sheaf $A$ on $\Aff /X$ is {\em weakly flasque} if, for any inclusion map $U\hookrightarrow V$ 
in $\Aff /X$ such that $U$ is a closed subscheme of $V$ defined by a nilpotent ideal, then 
the restriction $A(V)\rightarrow A(U)$ is surjective. Say that $E$ over $CRIS(X, \Ff )$ is
weakly flasque, if its restriction $E|_X$ is weakly flasque on $\Aff /X$. 

\begin{lemma}
\mylabel{fibrantweaklyflasque}
Suppose $A$ is a fibrant object in the injective model structure for complexes of sheaves of $\Oo$-modules over ${\Aff}/X$,
then each term of $A$ is weakly flasque. Also, locally free $\Oo$-modules are weakly flasque. 
\end{lemma}
\begin{proof}
Recall that fibrant objects are, in particular, complexes of injective sheaves, and injectivity
implies the required extension property. Similarly, locally free $\Oo$-modules satisfy the extension
property because we are working on the affine site. 
\end{proof}

Recall that the restriction of $\Diamond (\Lambda )$ to $\Aff /X$ is just $\Lambda$ considered as a quasicoherent sheaf on the affine site. 
Suppose $B$ is a sheaf of $\Oo$-modules on $CRIS(X, \Ff )$. Then we can form the complex of sheaves $\uHom^{\cdot} (\Diamond (\Lambda ),B)$
on $CRIS(X, \Ff )$. 

\begin{theorem}
\mylabel{uhompstar}
Suppose $B$ is a complex of
sheaves of $\Oo$-modules on the crystalline site $CRIS(X, \Ff )$ such that the restriction $B|_X$ is quasiisomorphic, as
a sheaf of $\Oo$-modules on ${\rm Aff}/X$, to a bounded complex of free $\Oo$-modules possibly of infinite rank. 
Suppose that $B$ is cohomologically flasque, and that the terms of $B$ are
weakly flasque. Then there is a natural quasiisomorphism 
of complexes of sheaves on $CRIS(X, \Ff )$
$$
\uHom^{\cdot} (\Diamond (\Lambda ),B)\stackrel{\sim}{\rightarrow} p_{\ast}(B|_X).
$$
\end{theorem}
\begin{proof}
We can localize on $X$ without affecting the conclusion. 
The map
$$
\uHom^{\cdot} (\Diamond (\Lambda ), B)|_X = \uHom^{\cdot} (\Diamond (\Lambda )|_X, B|_X) \rightarrow B|_X
$$
is given by evaluation at the unit section. By adjunction this gives the map in question. 

It suffices to prove the quasiisomorphism after restricting to $X$, because any object of $CRIS(X, \Ff )$ locally
admits a lifting into $X$. Put $A:= B|_X$, it is cohomologically flasque, weakly flasque and quasiisomorphic
to a bounded complex of locally free modules. 
Since $\uHom^{\cdot}$ commutes with restriction,
we are reduced to showing that the map
$$
\uHom^{\cdot} (\Diamond (\Lambda )|_X, A) \rightarrow (p_{\ast}(A))|_X
$$
is a quasiisomorphism.

These constructions have some invariance properties.  
If $A'\rightarrow A''$ is a quasiisomorphism of cohomologically flasque complexes of sheaves of $\Oo$-modules 
on ${\Aff}/X$ then it is objectwise a quasiisomorphism. 
On any open set where $\Diamond (\Lambda )|_X$ is an infinite rank free $\Oo$-module,
the internal $\uHom^{\cdot}$ from there to $A'$ becomes an infinite direct product of copies of $A'$.
Such a product preserves objectwise quasiisomorphisms so
$$
 \uHom ^{\cdot}(\Diamond (\Lambda )|_X, A')\stackrel{\sim}{\rightarrow}  \uHom ^{\cdot}(\Diamond (\Lambda )|_X, A'').
$$
Similarly, 
$$
(p_{\ast}(A'))(Y) = \lim _ {\leftarrow , k} A'(Y\times _X\Ff _k)
$$
where $\Ff _k=Spec ({\bf f}/{\bf j}^k)$ is the $k$-th formal neighborhood in $\Ff$. 
If $A'$ is weakly flasque then the maps in the inverse system are surjective. 
If $A'\rightarrow A''$ is a morphism between cohomologically flasque complexes of weakly flasque sheaves of $\Oo$-modules 
on ${\Aff}/X$ we get an inverse limit of objectwise quasiisomorphisms between inverse systems of surjective maps,
which is a quasiisomorphism (this is not the problematical situation detected in \cite{Neeman})
so $(p_{\ast}(A'))|_X \stackrel{\sim}{\rightarrow} (p_{\ast}(A''))|_X$.

Apply these starting with $A=B|_X$. Possibly after restricting to an open covering of $X$ we can choose a diagram
$A\rightarrow A' \leftarrow A''$
such that $A'$ is a fibrant object in the injective model structure for complexes of sheaves of $\Oo$-modules over ${\bf Aff}/X$,
and $A''$ is a bounded complex of locally free
$\Oo$-modules possibly of infinite rank. 
The latter conditions for $A''$ imply that
$$
\uHom^{\cdot} (\Diamond (\Lambda )|_X, A'') \rightarrow (p_{\ast}(A''))|_X
$$
is an isomorphism: expressing $\Diamond (\Lambda )= \bigcup _k\Diamond (J_k\Lambda )$ we get an inverse limit expression
on the left which is identical to the expression on the right. 

Now $A$ is cohomologically flasque and weakly flasque by hypothesis; and $A'$ and $A''$ also satisfy these conditions,
see Lemmas \ref{fibrantweaklyflasque} and \ref{cohflasque1}. 
The invariance statements described previously  imply that the map
$$
\uHom^{\cdot} (\Diamond (\Lambda )|_X, A) \rightarrow (p_{\ast}(A))|_X
$$
is a quasiisomorphism. 
\end{proof}

The above construction can be used to strictify complexes of sheaves on $CRIS(X, \Ff )$. 
In our notations for the affine case we use $\Lambda$ or $L=\Lambda (X)$ interchangeably:
$\Diamond (\Lambda )= \Diamond (L)$. The ring $L$ acts on the left on the complex of global sections
$$
\Hom ^{\cdot}(\Diamond (\Lambda ), B) = \Gamma (CRIS(X, \Ff ), \uHom^{\cdot} (\Diamond (\Lambda ), B))
$$ 
by the formula $\lambda \cdot f := f\circ (r_{\lambda})$
where the endomorphism $r_{\lambda}$  of $\Diamond (\Lambda )$ is right multiplication by $\lambda \in L$.
Hence
$\Hom ^{\cdot}(\Diamond (\Lambda ), B)$ becomes a complex of $L$-modules. In fact, this construction is exactly the
adjoint of the crystalization functor:
$$
\Diamond ^{\ast}(B) = \Hom ^{\cdot}(\Diamond (\Lambda ), B) = \Gamma ( CRIS(X, \Ff ), \uHom^{\cdot} (\Diamond (\Lambda ), B)).
$$

\begin{corollary}
Suppose $B$ is a fibrant complex of sheaves of $\Oo$-modules on $CRIS(X, \Ff )$ whose restriction to $\Aff /X$
is a perfect complex. Then there is a natural quasiisomorphism
of complexes of $A$-modules
$$
\Diamond ^{\ast}(B)
\stackrel{\sim}{\rightarrow} B(X).
$$
\end{corollary}
\begin{proof}
Theorem \ref{uhompstar} applies because
$B$ is cohomologically flasque by Lemma \ref{cohflasque1} and weakly flasque by Lemma \ref{fibrantweaklyflasque}. 
The global sections functor $\Gamma$ is the right adjoint of the
``constant sheaf'' functor which in turn preserves cofibrations and trivial cofibrations,
so they form a Quillen pair. It follows \cite{Quillen} \cite{Hovey} that $\Gamma$ preserves 
quasiisomorphisms between fibrant complexes.
The complex $\uHom^{\cdot} (\Diamond (\Lambda ), B)$ is fibrant.
On the other hand, $p_{\ast}$ is a right Quillen functor with left adjoint $p^{\ast}$, so 
$p_{\ast}$  preserves fibrant objects. By Theorem \ref{cohflares} $B|_X$ is fibrant, so 
$p_{\ast}(B|_X)$ is fibrant. 
Applying $\Gamma$ to the quasiisomorphism of Theorem \ref{uhompstar} gives
$$
\Gamma ( CRIS(X, \Ff ), \uHom^{\cdot} (\Diamond (\Lambda ), B))
\stackrel{\sim}{\rightarrow} \Gamma (CRIS(X, \Ff ), p_{\ast}(B|_X)) = B(X).
$$
\end{proof}

\begin{corollary}
\mylabel{Aperf}
In the affine case suppose $B$ is a fibrant complex of $\Oo$-modules on $CRIS(X, \Ff )$ whose pullback $B|_X$ is
globally quasiisomorphic to a strictly perfect complex.
Then $\Diamond ^{\ast}(B)$ is $\Oo _X$-perfect. 
\end{corollary}
\begin{proof}
The quasiisomorphism may be realized by a single map $P\rightarrow B|_X$.
Both $P$ and $B|_X$ are cohomologically flasque, so this induces a quasiisomorphism
$P(X)\stackrel{\sim}{\rightarrow} B(X)$. Thus $B(X)$ is $A$-perfect and the previous corollary gives
the conclusion.
\end{proof}

The global strictly perfect hypothesis in this corollary and the next theorem will be removed at the end of
the proof of \ref{diamondequiv} below.

\begin{theorem}
\mylabel{coherence1}
Suppose $X = \Spec  (A)$ is affine and $\Theta _{X,\Ff}$ is a free module. 
Suppose $B$ is a fibrant complex of $\Oo$-modules on $CRIS(X, \Ff )$ whose restriction to $\Aff /X$
is quasiisomorphic to a strictly perfect complex. 
Then there exists a bounded complex of finite rank locally free $L$-modules $R$ on 
$X$ and a quasiisomorphism $\Diamond (R)\stackrel{\sim}{\rightarrow} B$. 
\end{theorem}
\begin{proof}
By Corollary \ref{Aperf} and Proposition \ref{koszulcor} there is a resolution
$R\stackrel{\sim}{\rightarrow} \Diamond ^{\ast}(B)$
for $R$ a bounded complex of finite rank locally free $L$-modules. Note that \ref{koszulcor} works globally under
our assumptions on $B$ and $\Theta _{X,\Ff}$. By adjunction 
this corresponds to a map $\Diamond (R)\rightarrow B$ which we claim 
is a quasiisomorphism of complexes of $\Oo$-modules over $CRIS(X, \Ff )$.

Suppose $Y=\Spec  (C)\in CRIS(X,\Ff )$. We can lift (possibly after localizing on $Y$) to a map $Y\rightarrow X$. By definition 
$\Diamond (R)(Y) = C\otimes _AR$.
Let $P\rightarrow B|_X$ be the quasiisomorphism from a strictly perfect complex. Then $P(Y)=C\otimes _AP(X)$ and as in the proof of \ref{Aperf},
the map $P(X)\rightarrow B(X)$ is a quasiisomorphism of complexes of $A$-modules. 
We get quasiisomorphisms of perfect complexes of $A$-modules
$$
R = \Diamond (R)(X)  \rightarrow B(X) \leftarrow P(X).
$$

Choose a fibrant quasiisomorphism $M\rightarrow B(X)$ in the projective model category
$\cpx (A)$ of complexes of $A$-modules, with $M$ a cofibrant object. 
Note that both $R$ and $P(X)$ are bounded complexes
of projective $A$-modules, so we can lift to get a commutative diagram of the form
$$
\begin{array}{ccccc}
& & M & & \\
& \nearrow &\downarrow & \nwarrow \\
R  & \rightarrow & B(X) & \leftarrow & P(X) 
\end{array}
$$
where all maps are quasiisomorphisms. 

Next, notice that the functor $C\otimes _A -$ preserves quasiisomorphisms between cofibrant objects
in the projective model structure of $\cpx (A)$. 
Hence the upper diagonal arrows are quasiisomorphisms in the diagram
$$
\begin{array}{ccccc}
& & C\otimes _A M & & \\
& \nearrow & \downarrow & \nwarrow \\
C\otimes _A R  & \rightarrow & C\otimes _A B(X) & \leftarrow & C\otimes _A P(X) \\
\downarrow & &\downarrow & &\downarrow \\
\Diamond (R)(Y) & \rightarrow &B(Y) &\leftarrow &P(Y) 
\end{array}
$$
and also the lower squares commute, coming from the restriction maps from $X$ to $Y$
for the morphisms of complexes of sheaves $\Diamond (R)\rightarrow B \leftarrow P$. 

The left and right vertical arrows are equalities, because $P$ and $\Diamond (R)$ are strictly quasicoherent complexes. 

The same argument as before, saying that $P$ and $B|_X$ are both cohomologically flasque, implies that 
$P(Y)\rightarrow B(Y)$ is a quasiisomorphism. 

The right square and triangle imply that the map $C\otimes _AM\rightarrow B(Y)$ is a quasiisomorphism, then the
left square and triangle imply that the map $\Diamond (R)(Y)\rightarrow B(Y)$ is a quasiisomorphism. 
Note that the middle element of the diagram $C\otimes _AB(X)$ might not be quasiisomorphic to the rest, because
$B(X)$ itself may not be a cofibrant object for the projective model structure on $\cpx (X)$ (which was why we had to choose $M$ instead).

Putting these together for all $Y$ we get that the map $\Diamond (R)\rightarrow B$ is a quasiisomorphism. 
\end{proof}

\begin{theorem}
\mylabel{diamondequiv}
The functor $\Diamond $ induces an isomorphism on $Ext$ groups of bounded complexes of locally free $L$ modules: 
if $U$ and $V$ are bounded complexes of finite rank locally free $L$-modules then 
$$
Ext^i(\Diamond ): Ext ^i_{\mdl (L)}(U,V) \stackrel{\cong}{\rightarrow} Ext ^i_{\mdl (CRIS(X,\Ff ) , \Oo )}(\Diamond (U),\Diamond (V)).
$$
The extensions of $\Diamond $ to a DG-functor and simplicial functor respectively
$$
DG\cpx (L)\rightarrow DG\cpx (CRIS(X,\Ff ),\Oo )
$$
$$
L\cpx (L)\rightarrow L\cpx (CRIS(X,\Ff ),\Oo ) 
$$
are equivalences on the
subcategories of $\Oo$-perfect objects.
\end{theorem}
\begin{proof}
We can localize on $X$ using a spectral sequence. 
It suffices to prove the formula for $Ext$ when $U=V=L$ concentrated in degree $0$, and we can assume that $L$ is
a free $A$-module. Let $\Diamond (L )\rightarrow U'$ be a bounded below injective
resolution in $\mdl (CRIS(X,\Ff ) , \Oo )$, then $U'$ is fibrant in $\cpx (CRIS(X, \Ff )  , \Oo )$ and it satisfies
the hypotheses of Theorem \ref{uhompstar}, indeed $\Diamond(L)|_X$ is a free $\Oo$-module given by the free $A$-module $L$. 
Hence, 
$$
\Hom _{\cpx (CRIS(X, \Ff ) , \Oo )}(\Diamond (L ), U') \stackrel{\sim}{\rightarrow} U'(X)
$$
is a quasiisomorphism. On the other hand, the restriction of $U'$ to $\Aff /X$ is also a fibrant resolution of the restriction of
$\Diamond (L )$ (Theorem \ref{cohflares}), 
but $\Diamond (L )|_{\Aff /X}$ is the sheaf on the big (finite type, etale) site of $X$ associated to a 
quasicoherent sheaf. Since we are assuming that $X$ is affine, this has no higher cohomology so the resolution induces a quasiisomorphism
of sections over $X$:
$$
L = \Diamond (L )(X)\stackrel{\sim}{\rightarrow} U'(X)
$$
as desired. This proves that $\Diamond$ is fully faithful when restricted to the subcategory of bounded complexes of locally free $L$-modules. 

Applying the de Rham-Koszul resolution \ref{koszulcor}, any $\Oo_X$-perfect complex of $L$-modules is quasiisomorphic to a
bounded complexes of locally free $L$-modules possibly after localizing on $X$. So $\Diamond$ is Zariski locally
fully faithful on the full subcategory of $\Oo_X$-perfect objects. The $Hom$ objects on both sides are stacks
so it is fully faithful.  It takes an $\Oo _X$-perfect complex of $L$-modules to
an $\Oo $-perfect complex on $CRIS(X,\Ff )$. Between these categories its essential image contains anything whose pullback to $\Aff /X$ is globally quasiisomorphic
to a strictly perfect complex, by Theorem \ref{coherence1}. However, the source of the functor $\Diamond$ is a stack, and full faithfulness implies that its
image is a stack; therefore all $\Oo$-perfect complexes on $CRIS(X,\Ff )$ are in the essential image. 
\end{proof}


\section{Weak $\Lambda$-module structures on a given complex}
\mylabel{wpx}

\subsection{The partial Hochschild complex}

In this section we work in the affine case $X=\Spec (A)$, and $\Lambda$ corresponds to 
the filtered algebra $L = \Lambda (X)$. 
Use the tensor product $\otimes _A:= \otimes _{>A<}$ which contracts inner structures to define the algebra of tensor powers
$$
T^n L := L \otimes _{A} \cdots \otimes _{A} L \;\;\;\; TL := \bigoplus _{n\geq 0} T^nL,
$$
where $T^0L:= A$. These have both left and right $R$-module structures.
Maintain
the convention that $\otimes _A$ means $\otimes _{>A<}$ when that makes sense.

The tensor algebra has a partial Hochschild differential
$\delta : T^nL \rightarrow T^{n-1}L$ with $\delta ^2 =0$ defined by the formula
$$
\delta (u_1\otimes \cdots \otimes u_n):=
\sum _{i=1}^{n-1} (-1)^{i+1}u_1\otimes \cdots \otimes (u_iu_{i+1})
\otimes \cdots \otimes u_n.
$$
There is also a coassociative coproduct 
$$
\Psi :TL  \rightarrow TL \otimes _A TL 
$$
$$
\Psi (u_1\otimes \cdots \otimes u_n):=
\sum _{i=0}^{n} [u_1\otimes \cdots \otimes u_i]
\otimes
[u_{i+1} \otimes \cdots \otimes u_n]
$$
where by convention the empty bracket denotes 
the object $1\in T^0L$. 
The coproduct is compatible with the differential, using the appropriate signs. 

The filtration of $L$  induces a filtration 
$$
J_kTL = \bigoplus _n J_k(T^nL ), \;\;\;\; J_k(T^nL ):= \sum _{j_1+\ldots + j_n =k}
J_{j_1}L\otimes _A\cdots \otimes _A J_{j_n}L.
$$
The differential $\delta$ preserves the filtration, and coproduct is multiplicative on the filtration.

\subsection{The $Q(E,F)$ complex}
\mylabel{complexQ}
Suppose $(E,d_E)$ and $(F,d_F)$ are
complexes of $A$-modules. Define a complex of $A$-modules denoted $(Q(E,F), d_Q)$ as 
$$
Q(E,F):= \Hom _A(TL \otimes  _A E,F)
$$
with the differential $d_Q$ combining the differential $\delta$ on $TL$ 
with the differentials $d_E$ and $d_F$ using the usual sign rules. Note that the upper indexing for 
$TL$  with respect to which the differential is homological, is changed to negative indexing with
a cohomological differential when we do the construction of $Q$. Elements of $Q^i$ are multilinear functions
$$
q(u_1, \ldots , u_k; e)\in  F^{i+j-k} \;\;\; \mbox{for}\;\;\; u_l\in L, \;\; e\in E^j
$$
satisfying
$$
q(au_1, \ldots ; e) = aq(u_1 \cdots ; e),\;\;\; q(\ldots , u_ka; e) = q(\ldots , u_k; ae),
$$
$$
q(\ldots , u_ia , u_{i+1}, \ldots ; e) = q(\ldots , u_i , au_{i+1}, \ldots ; e) .
$$
The differential is 
\begin{eqnarray*}
(d_Q q)(u_1, \ldots ,u_k; e) &=& (-1)^{i+k}q(u_1, \ldots , u_k; d_Ee) \\
& +&  d_Fq(u_1, \ldots , u_k; e) \\
&+& \sum _{l=1}^{k-1}(-1)^{i+l+1}q(u_1, \ldots , (u_lu_{l+1}), \ldots , u_k; e).
\end{eqnarray*}
There is also $H: \Hom _A(T^{k+1}L\otimes _AE,F)\rightarrow \Hom _A(T^{k}L\otimes _AE,F)$ with
$$
(Hq)(u_1, \ldots , u_k; e):= \sum _{m=1}^{k}(-1)^{i+m}q(u_1, \ldots , \stackrel{m}{1} , \ldots , u_k; e),
$$
which anticommutes with $d_Q$.

Dual to the filtration of $TL$ let 
$J^kQ(E,F)$ consist of those morphisms $q$ which vanish on $J_{k-1}(TL)\otimes _AE$.
Note that $J^0Q(E,F)= Q(E,F)$ and $J^1Q(E,F)$ is the subcomplex
of elements $q$ which
restrict to zero on the first filtered piece of $TL$.
In view of our convention contracting the $A$-module
structures in the tensor algebra,  we have
$$
J_0(TL)=TA = T^0 A \oplus T^1A \oplus T^2A \oplus \cdots  = A\oplus A \oplus A \oplus \cdots .
$$
Thus, the quotient of $Q(E,F)$ by the subcomplex is
$$
Q(E,F)/J^1Q(E,F) \cong \Hom _A(TA\otimes _A E,F) .
$$
Denote the projection by
$$
P : Q(E,F)\rightarrow \Hom _A(TA\otimes _A E,F)
$$
and let $P_0: Q(E,F)\rightarrow \Hom _A(E,F)$ correspond to the piece $T^0A$. 

The coproduct $\Psi$ on $TL$ gives us an associative product
$$
Q(E,F)\times Q(F,G)\rightarrow Q(E,G),
$$
compatible with the differentials, denoted by simple juxtaposition in what follows.
It is unitary with respect to the identity morphisms of complexes, thought of as elements of
$Q(E,E)$. 

The filtration $J$ is preserved by the differential $d_Q$, and compatible with product. 
The product maps via $P_0$ to the usual composition on $\Hom _A(E,F)$. 
 
\subsection{Weak $L$-modules}
\mylabel{weakLmod}

Use the above constructions to define a differential graded category $\wpx (L)$ of 
{\em complexes with weak $L$-module structures}.  
 
Let $E$ be a complex of $A$-modules. The identity $1_E$ can
be thought of as an element of degree $1$ denoted
$$
1_E(1)\in \Hom _A(T^1A\otimes _AE, E) \subset \Hom _A(T^1A\otimes _AE, E).
$$

A {\em Maurer-Cartan element} (``MC-element'' for short) for a complex of $A$-modules $E$ is an element 
$$
\eta \in Q^1(E,E)
$$
such that $P(\eta )= 1_E(1)$, 
and such that $\eta$ satisfies the {\em MC equation}
$$
d_Q (\eta )+ \eta ^2 = 0.
$$ 
It is intended to represent a weak $L$-module structure on $E$.
The normalization condition $P(\eta )= 1_E(1)$ is there to insure that the action of $A=J_0L\subset L$
is equal to the given one. 

The condition on $P(\eta )$ implies that the coefficient of $\eta$ in $\Hom _A(T^0L\otimes _AE,E)$ is zero.
This corresponds to saying that $\eta$ doesn't deform the given differential of the complex $E$. Another possible notational 
choice would have been to integrate this differential into $\eta$ but we don't follow that route here. 

\begin{lemma}
\mylabel{weakify}
Suppose $\eta$ is an MC element for $E$. Then $\eta$ induces an action of $L$ on the
homology modules $H^i(E)$ compatible with the action of $A$.  
On the other hand, a strict left action of $L$ on $E$
(compatible with the given action of $A$) induces an MC element $\phi$ which gives back the
same action on $H^i(E)$. 
\end{lemma}
\begin{proof}
If $e$ represents a class $[e]\in H^i(E)$ and $u\in L$ then $\eta (u; e)$ represents the product
$u\cdot [e]\in H^i(E)$; the term $\eta (u,v; e)$ provides a coboundary showing that
$u\cdot (v\cdot [e]) = (uv)\cdot [e]$ thanks to the MC equation. Compatibility with the
action of $A$ comes from the condition $P(\eta )= 1_E(1)$. 

Suppose $E$ is a complex of $L$-modules. Then we can define an MC element by
$\phi (u;  e):= u\cdot e$ but $\phi (u_1, \ldots , u_k; e):= 0$ for $k\geq 2$. 
\end{proof}

If $\eta$ is a Maurer-Cartan element for $E$ and $\varphi$ is a Maurer-Cartan element for $F$,
then we define a new differential $d_{\eta , \varphi }$ on $Q(E,F)$ as follows:
$$
d_{\eta , \varphi }(q) := d_Q(q) + \varphi a - (-1)^{|q|}q\eta .
$$
One can verify using the MC equation that $d_{\eta , \varphi }^2=0$. 

An {\em MC complex} is a pair $(E,\eta )$ where $E$ is a complex of $A$-modules and $\eta$ is a 
Maurer-Cartan element. We will usually require that $E$ be bounded 
and consist of locally free $A$-modules in each degree.

Define the {\em weak morphism complex}
$$
\Hom _{\wpx (L)}((E,\eta ),(F,\varphi )):= \left( Q(E,F), d_{\eta , \varphi } \right) .
$$
There is a quotient map $P$ to  $\Hom _A(T^{\cdot}A\otimes _AE, F)$ and this can be projected
onto the piece $\Hom _A(T^0A\otimes _AE,F)$. This projection $P_0$ is compatible with the differential to give
a map of complexes
$$
P_0 : \left( \Hom _{\wpx (L)}((E,\eta ),(F,\varphi )) , d_{\eta , \varphi }\right) \rightarrow (\Hom _A(E,F), d_{E,F}).
$$
The image in $\Hom _A(E,F)$ of an element of $Q(E,F)$ will be called the {\em associated usual morphism}.

The {\em normalized weak morphism complex}  is the subcomplex 
$$
\Hom _{\wpx _0(L)}((E,\eta ),(F,\varphi ))\subset \Hom _{\wpx (L)}((E,\eta ),(F,\varphi ))
$$
consisting of elements $f$ such that $P(f)$ goes to zero in the higher 
terms of the form $\Hom _A(T^{m}A\otimes _AE, F)$, $m>0$. This is equivalent to requiring $P(f) = P_0(f)$. 

If $(E,\eta )$, $(F,\varphi )$ and $(G,\gamma )$ are
three MC complexes,  the composition map 
$$
Q(E,F)\otimes Q(F,G)\rightarrow Q(E,G)
$$
is still compatible with the new differentials 
$ d_{\eta , \varphi }$, $ d_{\varphi , \gamma }$ and $ d_{\eta , \gamma }$. 
We get a composition operation on weak morphism complexes (as well as their normalized subcomplexes)
which is compatible with the projections $P$ to associated usual morphisms. 

In what follows, the terminology ``locally free'' will include objects possibly of infinite rank. Note however that
there should be a single open covering on which the object becomes free, so an infinite direct sum of locally free objects is
not necessarily locally free. The de Rham-Koszul resolutions considered in \S \ref{koszul} are $A$-perfect bounded locally free complexes of $L$-modules. 

Let $DG\cpx ^{{\rm bplf}}(L)$
denote the differential graded category of bounded locally free complexes of $L$-modules which are $A$-perfect. 
Note that objects of $DG\cpx ^{{\rm bplf}}(L)$ are, in particular, complexes of locally free
$A$-modules. 
Similarly denote by $DG\cpx ^{{\rm bplf}}(A)$ the differential graded category of bounded complexes of
locally free $A$ modules which are perfect, and $\Pp ^{\rm str}$ the projection from $L$-modules to $A$-modules.

\begin{definition}
\mylabel{wpxdef}
Define $\wpx (L)$ to be the differential graded category of MC complexes $(E,\eta )$
with weak morphism complexes as above. Let $\wpx _0(L)$ denote the subcategory obtained by using the normalized morphism
complexes. 
Let  ${\wpx}^{{\rm bplf}}(L )$ and ${\wpx}^{{\rm bplf}}_0(L )$ be the full subcategories of objects
$(E,\eta )$ such that $E$ is bounded, locally free over $A$ and $A$-perfect.  Let 
$\Pp$ and $\Pp _0$ denote the projections from these categories to $DG\cpx ^{{\rm bplf}}(A)$ using 
$P$ on morphism complexes. The {\em weakification functor} 
$$
\Ww _0: DG\cpx ^{{\rm bplf}}(L)  \rightarrow{\wpx}^{{\rm bplf}}_0(L ) 
$$
represents the operation of the second statement in Lemma \ref{weakify} and $\Ww$ denotes its composition 
with the inclusion of the normalized subcategory. 
\end{definition}

\begin{lemma}
\mylabel{normequiv}
The inclusion $\nu : {\wpx}^{{\rm bplf}}_0(L )\hookrightarrow {\wpx}^{{\rm bplf}}(L )$ is a quasiequivalence. 
\end{lemma}
\begin{proof}
The complex $\bigoplus _{m>0}\Hom _A(T^{m}A\otimes _AE, F)$ with differential determined by $1_E(1)$ and $1_F(1)$ is acyclic: its differentials are alternately
isomorphisms or zero. This is the quotient by inclusion of normalized into non-normalized morphism complexes. 
\end{proof}

\begin{lemma}
\mylabel{qisinvariance}
A closed degree zero morphism $(E,\eta )\stackrel{f}{\rightarrow}(F,\varphi )$ in ${\wpx}^{{\rm bplf}}(L )$
is an inner equivalence if and only if the underlying usual morphism
of complexes $P_0(f)$ is a quasiisomorphism. 
\end{lemma}
\begin{proof}
If $P_0(f)$ is a quasiisomorphism then for any MC-complex
$(G,\gamma )$ in ${\wpx}^{{\rm bplf}}(L )$, composition with $f$ induces quasiisomorphisms (with respect to the MC-twisted differentials)
$$
Q(F,G)\stackrel{\sim}{\rightarrow} Q(E,G), \;\;\; 
Q(G,E)\stackrel{\sim}{\rightarrow} Q(G,F).
$$
This can be seen by putting the trivial filtrations on $E$, $F$ and $G$ and passing to the associated-graded objects. 
In the other direction if $f$ is an inner equivalence then so is $P_0(f)$. 
\end{proof}

The following homotopy invariance gives the fibrant property for the functor $\Pp _0$ which explains why we are interested in
the notion of weak structure.

\begin{lemma}
\mylabel{maininvariance}
Suppose $(E,\eta )$ is an MC complex, and suppose $F$ is a complex of $A$-modules. Suppose
both are bounded complexes of locally free $A$-modules. Suppose 
$a_0:E\rightarrow F$ is a quasiisomorphism of complexes of $A$-modules.
Then there exists an MC element $\varphi$ for $F$, and 
$$
a\in Q^0(E,F), \;\;\; P(a)=P_0(a) = a_0, \;\;\; d_{\eta ,\varphi}(a)=0.
$$
The same for an equivalence going in the other direction. 
\end{lemma}
{\em Proof:}
By the condition that the complexes are bounded and locally free over $A$, we can choose a homotopy inverse 
$b_0:F\rightarrow E$ to $a_0$, with $K\in End (E)^{-1}=\Hom ^{-1}(T^0A\otimes _AE,E)\subset Q(E,E)^{-1}$ 
with $d(K) = b_0a_0 - 1$.
Put 
$$
\varphi := a_0 \eta (\sum _{m=0}^{\infty}(K\eta ) ^m) b_0.
$$
The sum converges in $Q(F,F)$, indeed $K\eta (;e)=0$ so $(K\eta )^m$ vanishes on $T^kL\otimes E$ for $k<m$,
and for $m\gg 0$ such a thing which is also in $Q(E,E)^0$ must vanish since $E$ is a bounded complex. 
We have 
$$
d_Q(\varphi ) = \sum \pm a_0 \eta \cdots Kd(\eta )K\ldots \eta b_0 + 
\sum \pm a_0 \eta \cdots \eta d(K)\eta \ldots \eta b_0 . 
$$ 
It is left to the reader to fill in the signs.
Using $d(K) = b_0a_0-1$ and $d_Q(\eta )=-\eta ^2$, the terms corresponding to
$b_0a_0$ in $d(K)$ give $-\varphi ^2$ whereas the terms corresponding to $-1$ in 
$d(K)$ and to $d(\eta ) = -\eta ^2$ cancel. Thus,
$$
d_Q(\varphi ) = -\varphi ^2.
$$
Similarly,  put
$$
a := a_0 (\sum _{m=0}^{\infty}(\eta K) ^m).
$$
The inner terms in $d(a)$ work the same way as before; the term involving $d(K)$ on the
right end gives $a\eta$ so we get
$$
d(a) +\varphi a - a \eta = 0.
$$
This says that $a$ is a morphism from $\eta$ to $\varphi$.
\eop

\begin{corollary}
\mylabel{Pfibration}
The functor $\Pp _0: {\wpx}^{{\rm bplf}}_0(L )\rightarrow DG\cpx ^{{\rm bplf}}(A)$
is a fibration of dg categories in the sense of \S \ref{htyfib}. 
\end{corollary}
\begin{proof}
It is clearly surjective on mapping complexes, and the lifting condition is the previous lemma.
\end{proof}

Our main strictification result says that weak complexes are equivalent to strict complexes. 

\begin{theorem}
\mylabel{upsiloncoherence}
The functors $\Ww_0$ or $\Ww$ are quasi-equivalences of DG categories relative to 
$DG\cpx ^{{\rm bplf}}(A)$ via the projections $\Pp ^{\rm str}$ and $\Pp _0$ or $\Pp$.
\end{theorem}

In view of Lemma \ref{normequiv} we can use $\Ww$ although the application will concern $\Ww _0$.
We use the strategy suggested by Lemma \ref{genadjoint}, by defining a right adjoint $\Ww ^{\ast}$ to $\Ww$.
In fact we are in pretty much the same situation as before: even though $\Ww$ doesn't 
come from a functor on abelian module categories, it still starts
from a category of complexes of modules over $L$. 

Recall  that $\Ww (L) = (L,\phi )$ where $\phi$ is the MC element for $L$
constructed in Lemma \ref{weakify}, and $\Ww (L)$ admits
$L$ as an algebra of right endomorphisms. For an MC-complex $(E,\eta )$ define
$$
\Ww ^{\ast}(E,\eta ) := \Hom _{\wpx (L)}(\Ww (L), (E,\eta ))= \left( Q(L,E) , d_{\phi , \eta} \right) 
$$
with left action of $L$. It is the simple complex associated to the double complex 
$$
\Hom _A(T^1L, E) \rightarrow \Hom _A(T^2L,E)\rightarrow \ldots .
$$
Complete this  by adding in $E$ at the beginning
$$
E\rightarrow \Hom _A(T^1L, E) \rightarrow \Hom _A(T^2L,E)\rightarrow \ldots 
$$
and denote the resulting full complex by $(Q(L,E)^+, d^+_{\phi , \eta})$ 
for the natural extension of the differential $d_{\phi , \eta}$.

\begin{lemma}
\mylabel{wkLacyclic}
For any $(E,\eta )$, the augmented complex $(Q(L,E)^+, d_{\phi , \eta})$ is acyclic. This gives a 
quasiisomorphism of complexes of $A$-modules 
$$
E\stackrel{\sim}{\rightarrow} \left( Q(L,E) , d_{\phi , \eta} \right) .
$$
\end{lemma}
\begin{proof}
The elements of $Q(L,E)^+$ are functions with the usual multilinearity properties written $q(u_1,\ldots , u_k)$. 
Use an extended version $H^+$ of the homotopy defined above which inserts $1$ at each place. 
Calculate $(d^+_{\phi , \eta}H^+ + H^+d^+_{\phi , \eta})q = q\cdot H(\eta )$. However, $H(\eta )$ is an element of $Q(E,E)$ which
starts off with the identity, due to the condition $P(\eta )= 1_E(1)$. The remainder $H(\eta )- 1$ involves only positive degree tensors, so it is nilpotent
by boundedness of $E$. Thus $H(\eta )$ is an invertible endomorphism homotopic to zero, so $(Q(L,E)^+, d_{\phi , \eta})$ is acyclic.
\end{proof}

The following corollary shows that our functor $\Ww$ is quasi-fully-faithful. 

\begin{corollary}
\mylabel{mapinvariance}
Suppose $E$ and $F$ are bounded complexes of $L$-modules with corresponding MC elements $\eta$ and $\varphi$ respectively, 
such that in each degree 
$E$ is a locally free $L$-module. Then the map
$$
\Hom ^{\cdot}_{L}(E,F)\rightarrow (Q(E,F), d_{\eta ,\varphi })
$$
is a quasiisomorphism. 
\end{corollary}
{\em Proof:}
Both sides satisfy a homotopy Zariski localization property (for the right hand side, this can be shown using the technique of
\S \ref{sub-acyclic} below), so
it suffices to consider the case $E\cong L$. Then $T^iL\otimes _AL = T^{i+1}L$.  The map of complexes in question becomes
$$
F\rightarrow \left( Hom _A(T^1L , F) \rightarrow Hom _A(T^2L , F) \rightarrow \ldots \right) .
$$
The total complex including $F$ is the same as that of Lemma \ref{wkLacyclic}, so it is acyclic. 
This says that our map of the present lemma is a quasiisomorphism. 
\eop

\begin{corollary}
\mylabel{Qperfect}
If $E$ is perfect as a complex of $A$-modules, then $\Ww ^{\ast}(E,\eta )= (Q(L,E) , d_{\phi , \eta})$ is perfect as a complex of $A$-modules.
\end{corollary}
\begin{proof}
By Lemma \ref{wkLacyclic}, they are quasiisomorphic. 
\end{proof}

We now prove essential surjectivity of $\Ww$.
Note  that $(Q(L,E) , d_{\phi , \eta})$ is in fact a complex of $L$-modules: the right multiplication by $L$ on $T^iL$ (for $i\geq 1$) transforms to a left
multiplication on $\Hom _A(T^iL, E)$. 
The de Rham-Koszul resolution of Proposition \ref{koszulcor} now applies: we can choose a resolution 
$$
R\stackrel{\sim}{\rightarrow} \left( Q(L,E) , d_{\phi , \eta} \right)
$$
by a bounded complex of finite rank locally free $L$-modules. 

As part of the ``adjunction'' between $\Ww$ and $\Ww ^{\ast}$, there is a tautological closed degree zero element of the weak morphism complex
$$
\xi \in  \Hom ^0_{\wpx (L)}( \Ww ( Q(L,E) , d_{\phi , \eta}),(E,\eta ))= Q(Q(L,E),E)^0
$$
defined by the formula
$\xi (u_1, \ldots , u_i ; f):= f(u_{1},\ldots , u_i; 1)$.
The differential of $\Hom _{\wpx (L)}(\Ww ( Q(L,E) , d_{\phi , \eta}),(E,\eta ))$ is $d_{\psi , \eta}$ 
where $\psi$ denotes the MC structure from Lemma \ref{weakify} for the $L$-module $Q(L,E)$. 
Calculate
$$
d_{\psi , \eta}(\xi ) = d_Q\langle d_{\phi , \eta}, d_E\rangle (\xi ) + 
\eta \xi - (-1)^{|\alpha |}\xi \psi = 0.
$$
Furthermore the underlying morphism of complexes of $A$-modules corresponding to $\xi$ is a quasiisomorphism, indeed it is 
a strict left inverse to the quasiisomorphism in Lemma \ref{wkLacyclic} above. 
We can then compose $\xi$ with the image by $\Ww$ of the morphism $R\rightarrow ( Q(L,E) , d_{\phi , \eta} )$
to get an MC-morphism inducing a usual quasiisomorphism of $A$-modules
$\Ww (R) \stackrel{\sim}{\rightarrow} (E,\eta )$.
This quasiisomorphism is an inner equivalence in the dg-category $\wpx ^{\rm bplf}(L)$, so it completes the proof of 
essential surjectivity to finish the proof of Theorem \ref{upsiloncoherence}.
\eop


\section{\v{C}ech globalization}
\mylabel{cechglob}

\subsection{Simplicial globalization of the crystalization functor}
\mylabel{sec-lambdamodel}

We now extend Theorem \ref{diamondequiv} to a formal category $(X,\Ff )$ where $X$ is not affine. 
For any etale map $Y\rightarrow X$, the pullback of $\Lambda$ to $Y$ is again a split almost-polynomial 
sheaf of rings of differential operators $\Lambda _Y$ associated to a formal category of smooth type $(Y,\Ff _Y)$ where
$$
\Ff _Y := \left( \Ff \times _{X\times X}Y\times Y \right) _{\rm diag}.
$$
The subscript ${\rm diag}$ means that we take the open and closed sub-formal scheme supported along the diagonal of $Y\times Y$.

\begin{lemma}
\mylabel{lambdastack}
The $(\infty , 1)$-prestacks on the small etale site $X_{\rm et}$ defined by
$$
\left[ Y\mapsto  L\cpx _{\rm qc}(\Lambda _Y)\right] \;\;\;\mbox{and} \;\;\; 
\left[ Y\mapsto  L\cpx (CRIS(Y,\Ff _Y ),\Oo  ) \right]
$$
are stacks.  The crystallization functors $\Diamond _Y$  defined for affine $Y$ piece together
by this stack condition
to give a crystallization functor 
$$
\Diamond : L\cpx _{\rm qc}(\Lambda )\rightarrow L\cpx (CRIS(X,\Ff ),\Oo ).
$$
Restricting to the full subcategories of $\Oo$-perfect objects on both sides, $\Diamond$ becomes an equivalence. 
\end{lemma}
\begin{proof}
Use the criterion \cite[Theorem 19.4]{HirschowitzSimpson} on the left Quillen presheaves
$Y\mapsto \cpx _{\rm qc}^{\rm inj}(\Lambda _Y)$ and $Y\mapsto \cpx _{\rm qc}^{\rm inj}(CRIS(X,\Ff ),\Oo )$
over the small etale site $X_{\rm et}$.

For affine $Y$  the crystallization functor $\Diamond _Y$ was defined by first going to the projective model category.
These satisfy naturality with respect to restrictions. The stack condition gives in particular homotopical descent for 
the category $\Uu$ of all affine Zariski open sets of $X$, so  
$L\cpx _{\rm qc}(\Lambda )$ is the homotopy limit over $\Uu$ of $\left[ Y\mapsto  L\cpx _{\rm qc}(\Lambda _Y)\right] $
and similarly $L\cpx (CRIS(X,\Ff ),\Oo )$ is the homotopy limit over $\Uu$ of $\left[ Y\mapsto  L\cpx (CRIS(Y,\Ff _Y ),\Oo  ) \right]$.
The limit of functors $\Diamond _Y$ gives $\Diamond$.

The restrictions of $\Diamond _Y$ to subcategories of $\Oo$-perfect objects are equivalences by Theorem \ref{diamondequiv}.
Thus the same holds for $\Diamond$. 
\end{proof}

Suppose given a morphism of smooth type $f$ from $(X,\Ff )$ to a base scheme $S$. If $Z\rightarrow S$ is a morphism of schemes of finite type, 
let $X_Z:=X\times _S Z\rightarrow Z$ and the pullback of $\Lambda$ to $X_Z$ is the base change $\Lambda _Z$ 
defined using the fact that $f^{-1}(\Oo _S)$ is in the center of $\Lambda$ \cite[Lemma 2.]{moduli1}.
We would like to obtain a presheaf of simplicial categories
$$
 Z\mapsto L\cpx _{\Oo -{\rm perf}}(\Lambda _Z)\cong \widetilde{DP}DG\cpx _{\Oo -{\rm perf}}(\Lambda _Z) 
$$
on the category of affine schemes of finite type $Z\in \Aff / S$. As before, this will be problematic if we use the injective
model category structure, since $Z$ may not be flat over $S$. To calculate pullbacks  for
maps $Z'\rightarrow Z$ in $\Aff /S$ we should use the projective model structure. Hence we need to assume that $X\rightarrow S$
is an affine morphism. These can then be patched together using a stack condition.

\begin{lemma}
\mylabel{functoriality}
If $X\rightarrow S$ is an affine morphism, we obtain a presheaf of simplicial categories $[Z\mapsto L\cpx (\Lambda _Z)]$ 
which is a stack on the site $\Aff /S$. If $X\rightarrow S$ is quasi-separated of finite type and flat,
one can cover it by open subsets which are affine over $S$ and take the inverse limit of the resulting system of  stacks
to obtain a stack on $\Aff /S$ again denoted $[Z\mapsto L\cpx (\Lambda _Z)]$. 
\end{lemma}
\begin{proof}
For the stack condition it 
suffices to have descent for an etale covering $Z'\rightarrow Z$. In this case, $X\times _SZ'\rightarrow X\times _SZ$ is
an etale covering in the small etale site over $X\times _SZ$, so the previous descent Lemma \ref{lambdastack} applies. 
\end{proof}

\subsection{DG globalization of complexes}
\mylabel{dgglobcpx}

There is a multiplicative \v{C}ech resolution for a sheaf of differential graded algebras. 

Suppose $X$ is a topoligical space, and $\Uu$ is an open covering. Think of $\Uu$ as being the semicategory of
multiple intersections in the covering, where the morphisms are only the nontrivial inclusions
(throw out the identity inclusions). Note that there is at most one morphism between elements of $\Uu$,
so we leave the morphisms out of our notation.

Suppose $A$ is a complex of sheaves on $X$, or even just a presheaf of complexes over the category $\Uu$. 

Define the {\em local sections of $A$ over $\Uu$} to be the following complex of groups denoted $G=G_{\Uu}A$.
Let $G^i$ be the set of functions 
$$
g (U_0,U_1,\ldots ,U_k)\in A^{i-k}(U_0)
$$
defined whenever $U_0 \subset U_1 \subset \ldots \subset U_k$ is a strictly increasing sequence of objects of $\Uu$.
Set
$$
(dg) (U_0,\ldots , U_k):= d (g (U_0,\ldots , U_k)) + \sum _{j=0}^k (-1)^j g(U_0, \ldots , \widehat{U_j},
\ldots , U_k)|_{U_0}.
$$
The restriction to $U_0$ is necessary only for the single term $j=0$ where the value starts out in $A(U_1)$ and needs to be 
restricted to $U_0$. 
This defines a differential with $d^2=0$. 

We have a product 
$$
\mu : G_{\Uu}(A)\otimes G_{\Uu} (B) \rightarrow G_{\Uu}(A\otimes B)
$$
defined by 
$$
\mu (f\otimes g) (U_0,\ldots , U_k) := \sum _{j=0}^k f(U_0,\ldots , U_j)\otimes g(U_j,\ldots , U_k)|_{U_0}.
$$
This is associative, and compatible with the differential (for the same reason as before). In particular, if
$A$ is a presheaf of differential graded algebras, then $G_{\Uu}(A)$ has a natural structure of differential
graded algebra.

\subsection{A sheaf-theoretic version}
\mylabel{sheafglob}

We can also define a sheaf-theoretic version of this construction, obtained by replacing $A^{i-k}(U_0)$
by the direct image $j_{U_0/X}(A^{i-k}|_{U_0})$ in the above definition. Call this $\Gg _{\Uu}(A)$.
A section $a$ of $A^i$ gives a section of $\Gg _{\Uu}(A)$ obtained by setting 
$g(U_0):= a$ and $g(U_0,\ldots , U_k):= 0$ for $k>0$. This is compatible with the differential, so it gives a map
of complexes of presheaves 
$$
A\rightarrow \Gg _{\Uu}(A).
$$ 
It is easy to verify that this is a quasiisomorphism. Indeed, if $X$ appears as part of the covering $\Uu$ then 
$A(X)\rightarrow G_{\Uu}(A)$ is a quasiisomorphism by a classical calculation. Therefore in general the above map of
presheaves of complexes induces a quasiisomorphism on any open subset contained in some element of the covering. 

Suppose $X$ is a quasi-separated scheme and $\Uu$ is an affine open covering (in particular all of the open sets in $\Uu$ 
which includes the multiple intersections, are affine). Suppose that $A$ is a complex of quasicoherent sheaves.
Then the elements of $\Gg _{\Uu}(A)$ are direct images of quasicoherent sheaves via affine inclusions, so they are
acyclic. In particular, 
$$
A\rightarrow \Gg _{\Uu}(A)
$$ 
is an acyclic resolution. Again if $A$ is a sheaf of quasicoherent differential graded algebras then 
$\Gg _{\Uu}(A)$ is a quasicoherent differential graded algebra whose components are acyclic,
and the map is a map of sheaves of dga's. 

Note that $G_{\Uu}(A)$ is the complex (or dga) of global sections of $\Gg _{\Uu}(A)$.

It is instructive to write down this resolution for the example of an open covering with two elements 
denoted $U,V$, for the complex $A:= \Oo$.  Denote by $UV$ the intersection. Denote for example by $\Oo_{UV}$ the
direct image from $UV$ to $X$ of the sheaf $\Oo$.
Then our resolution $\Gg _{\Uu}(\Oo )$ takes the form
$$
\Oo_U \oplus \Oo_V \oplus \Oo_{UV} \rightarrow  \Oo_{UV} \oplus \Oo_{UV}.
$$
It is a little bit bigger than the standard \v{C}ech resolution, but equivalent,
the difference being the acyclic complex $\Oo _{UV} \rightarrow \Oo _{UV}$. 

\subsection{Globalization of dgc's}
\mylabel{dgcglob}

In a similar way we will define a \v{C}ech globalization of a presheaf of differential graded categories.
Suppose $C$ is a presheaf of dgc's over the category $\Uu$. 
Then define a differential graded category $G_{\Uu}(C)$ as follows.
An object is a pair $(E,\eta )$ where $E$ is a collection of objects $E(U)\in {\rm ob}(C (U))$,
and $\eta$ a collection of elements $\eta (U_0, \ldots , U_k)$ of $Hom ^{1-k}_{C(U_0)}(E(U_0), E(U_k)|_{U_0})$
indexed by strictly increasing sequences $U_0\subset \ldots \subset U_k$ in $\Uu$ (for $k\geq 1$). This is 
subject to a Maurer-Cartan equation of the form $d(\eta ) + \eta ^ 2=0$, where the differential and product are
defined much as previously. Define the complex of morphisms
$$
Hom _{G_{\Uu}(C)}((E,\eta ) , (F,\varphi ))
$$
to be the complex whose piece of degree $i$ consists of collections of functions $a(U_0,\ldots , U_k)
\in Hom ^{i-k}(E(U_0), F(U_k)|_{U_0})$,
with differential denoted $d_{\eta , \varphi}$ obtained by a formula analogous to the previous ones,
and with composition product defined as before too. 

Let $G_{\Uu}^{{\rm eq}}(C)$ denote the subcategory of objects where the principal restriction maps are equivalences, i.e. the
$\eta (U_0,U_1)$ are equivalences from $E(U_0)$ to $E(U_1)|_{U_0}$  in the dgc $C(U_0)$.

If $(C,\varepsilon )$ is a presheaf of augmented dgc's then define the {\em augmented globalization}
$G_{\Uu }(C,\varepsilon )$ to be the subcategory of objects of $G_{\Uu} (C)$ such that the transition maps
are mapped to $1$ by the augmentation; and with morphisms being those which map to a constant in $\cc$ by the augmentation. 
In our applications this condition will automatically put us into $G_{\Uu}^{{\rm eq}}(C)$.

\begin{lemma}
\mylabel{mcglobcommute}
The augmented globalization $G_{\Uu }(U\mapsto {\bf MC}(Q(U),\varepsilon ),\varepsilon )$ is equal 
to the differential graded category of
Maurer-Cartan elements of the globalization ${\bf MC}(G_{\Uu }(Q),\varepsilon )$.
\end{lemma}
\eop

\begin{lemma}
\mylabel{gucinvariant}
If $C\rightarrow C'$ is a morphism of presheaves of augmented dgc's over $\Uu$
which induces a quasiequivalence over each
element of $\Uu$, then $G_{\Uu }(C ,\varepsilon )\rightarrow G_{\Uu}(C',\varepsilon )$ is a quasiequivalence.
\end{lemma}
{\em Proof:} 
Since the covering is finite, there is a finite amount of data of the form $\eta (U_0,\ldots , U_k)$ to consider.
Also we only consider strict inclusions of open sets so the multiplication (even of the degree zero piece) in 
something like the formula $d(\eta ) + \eta ^ 2$, is nilpotent. Using this one can obtain the invariance.
\eop

{\em Remark:}  Clearly one can define $\Gg _{\Uu}(C,\varepsilon )$ when $C$ is defined only for the elements
of the covering $\Uu$, for example $C$ might only be defined for affine Zariski-open sets of a quasi-separated scheme.
The invariance result of Lemma \ref{gucinvariant} holds also in this case.

Let $[U\mapsto DG\cpx ^{\rm inj}_{\Oo -{\rm perf}}(\Lambda |_U )]$ denote the presheaf of d.g.c.'s on the Zariski topology of $X$, which 
associates to $U\subset X$ the d.g.c. of complexes of $\Lambda |_U$-modules which are cofibrant and fibrant objects in the injective model
structure, and which are $\Oo$-perfect.

\begin{lemma}
\mylabel{toledotong}
For any affine
open covering $\Uu$, the map
$$
DG\cpx ^{\rm inj}_{\Oo -{\rm perf}}(\Lambda _X) \rightarrow 
G _{\Uu}^{\rm eq}[U\mapsto DG\cpx ^{\rm inj}_{\Oo -{\rm perf}}(\Lambda |_U )]
$$
is a quasi-equivalence of differential graded categories. 
\end{lemma}
{\em Proof:}
Given a \v{C}ech-twisted complex, one can define in a natural way  its complex of sections over an open set. This 
has a structure of complex of $\Lambda$-modules, and the natural map from the original object to the new one is
seen to be a quasiisomorphism over any open subset contained in some element of the covering (according to the
usual principle that \v{C}ech-type resolutions including the full space are acyclic). 
\eop

This result fits in with the fact that $[U\mapsto \widetilde{DP}(DG\cpx ^{\rm inj}_{\Oo -{\rm perf}}(\Lambda |_U ))]$ is
a stack.

On the other hand, let 
$[U\mapsto DG\cpx ^{\rm bplf}(\Lambda |_U)]$ denote the presheaf of d.g.c.'s over the subcategory of affine open sets in $X$,
which to an affine $U$ associates the d.g.c. whose objects are bounded complexes of
locally free $\Lambda |_U$-modules which are $\Oo _U$-perfect.

For any affine $U$ let
$DG\cpx ^{{\rm bplf}/{\rm inj}}(\Lambda |_U)$ denote the dgc whose objects are morphisms
$$
E\rightarrow E'
$$
where $E$ is a bounded complex of locally free $\Lambda |_U$-modules (in particular, a cofibrant and fibrant object in
the projective model structure) and $E'$ is a fibrant (and automatically cofibrant) object in the injective
model structure; and the morphism between them is a quasiisomorphism of complexes of $\Lambda|_U$-modules. 
The structure of dgc is defined as follows:
$$
\Hom ^{\cdot}(E\rightarrow E', F\rightarrow F'):= \Hom ^{\cdot}(E,F)\oplus \Hom ^{\cdot}(E',F')\oplus \Hom ^{\cdot}(E,F')[1]
$$
with differential given by the mapping cone formula and composition defined by 
$(a,b,c)\circ (u,v,w) := (au, bv,bw + cu)$.  We have quasiequivalences of d.g.c.'s
\begin{equation}
\label{twoway}
DG\cpx ^{\rm bplf}(\Lambda |_U) \stackrel{\sim}{\leftarrow} 
DG\cpx ^{{\rm bplf}/{\rm inj}}(\Lambda |_U) \stackrel{\sim}{\rightarrow}
DG\cpx ^{\rm inj}_{\Oo -{\rm perf}}(\Lambda |_U )
\end{equation}
varying naturally with $U$. 
These induce quasi-equivalences on the globalizations, by Lemma \ref{gucinvariant}.
We obtain a diagram of quasi-equivalences relating $DG\cpx ^{\rm inj}_{\Oo -{\rm perf}}(\Lambda _X) $
to $G _{\Uu}^{\rm eq}[U\mapsto DG\cpx ^{\rm bplf}(\Lambda |_U)]$.

\subsection{Globalization of weak complexes}
\mylabel{globwkcpx}

For any $U\in \Uu$,  we have the  normalized dg category of weak complexes 
$\wpx ^{\rm bplf}_0(\Lambda (U))$ considered above
(Definition \ref{wpxdef}). By Theorem \ref{upsiloncoherence} the functor
$$
DG\cpx ^{\rm lf}_{\Oo -{\rm perf}}(\Lambda |_U) \rightarrow \wpx ^{\rm bplf}_0(\Lambda (U))
$$
is a quasiequivalence. Lemma \ref{gucinvariant} gives a quasiequivalence of globalizations
$$
G _{\Uu}^{\rm eq}[U\mapsto DG\cpx ^{\rm lf}_{\Oo -{\rm perf}}(\Lambda |_U)] \stackrel{\sim}{\rightarrow}
G _{\Uu}^{\rm eq}[U\mapsto \wpx ^{\rm bplf}_0(\Lambda (U)] .
$$
Combined with the previous diagram and Lemma \ref{toledotong}
we obtain a diagram of quasi-equivalences written vertically, and compatible with the corresponding diagram for
$\Oo$-modules:
$$
\begin{array}{ccc}
DG\cpx ^{\rm inj}_{\Oo -{\rm perf}}(\Lambda _X) 
& \rightarrow & DG\cpx ^{\rm inj}_{\Oo -{\rm perf}}(\Oo _X) \\
\downarrow & &\downarrow \\
G _{\Uu}^{\rm eq}[U\mapsto DG\cpx ^{\rm inj}_{\Oo -{\rm perf}}(\Lambda |_U)]
& \rightarrow &G _{\Uu}^{\rm eq}[U\mapsto DG\cpx ^{\rm inj}_{\Oo -{\rm perf}}(\Oo _U)] \\
\uparrow & &\uparrow \\
G _{\Uu}^{\rm eq}[U\mapsto DG\cpx ^{{\rm bplf}/{\rm inj}}_{\Oo -{\rm perf}}(\Lambda |_U)]
& \rightarrow &G _{\Uu}^{\rm eq}[U\mapsto DG\cpx ^{{\rm bplf}/{\rm inj}}_{\Oo -{\rm perf}}(\Oo  _U)] \\
\downarrow & &\downarrow \\
G _{\Uu}^{\rm eq}[U\mapsto \wpx ^{\rm bplf}_0(\Lambda (U))] 
& \rightarrow &G _{\Uu}^{\rm eq}[U\mapsto DG\cpx ^{\rm bplf}(\Oo (U))] .
\end{array} 
$$

\subsection{Description of the homotopy fiber}
The \v{C}ech globalization commutes with the dgc fiber when we have a collection of fibrant functors
of dgc's over the open sets of the covering.

\begin{lemma}
\mylabel{globfib}
If $A\rightarrow B$ is a morphism of presheaves of dgc's which is fibrant over each open set of 
$\Uu$, then $G_{\Uu}(A)\rightarrow G_{\Uu}(B)$ is fibrant. For a global section $E$ of $B(X)$,
the fiber and globalization operations commute.
\end{lemma}
\eop

This applies to the bottom horizontal map in the above big diagram. If $E$ is a perfect complex of $\Oo _X$-modules,
we can realize it globally as a bounded complex of locally free $\Oo_X$-modules \cite{ThomasonTrobaugh}, SGA VI. It gives elements of all of the different
\v{C}ech globalizations on the right in the big diagram. Hence, the augmented dg homotopy fiber of the upper map 
$$
DG\cpx ^{\rm inj}_{\Oo -{\rm perf}}(\Lambda _X) 
\rightarrow  DG\cpx ^{\rm inj}_{\Oo -{\rm perf}}(\Oo _X)
$$
over $E$, is equivalent to the strict fiber $\Fib ^+$ of the bottom map
\begin{equation}
\mylabel{mainmap}
G _{\Uu}^{\rm eq}[U\mapsto \wpx ^{\rm bplf}_0(\Lambda (U))] 
\rightarrow G _{\Uu}^{\rm eq}[U\mapsto DG\cpx ^{\rm bplf}(\Oo (U))] 
\end{equation}
over $E$. 
The explicit description of the fiber of \eqref{mainmap}
is the basic tool we use in the next section to check geometricity.


\section{Geometricity}
\mylabel{geomsec}

\subsection{Geometric stacks}
\mylabel{geomstacks}

The notion of {\em geometric $n$-stack} \cite{geometricN} is a direct generalization of the notion of Artin $1$-stack \cite{Artin}.
It has been generalized to a notion of geometric derived stack in \cite{hag2}. For our purposes a {\em locally geometric $(\infty ,0)$-stack}
will mean an $(\infty ,0)$-stack  
which is locally a geometric $n$-stack of groupoids for some $n$ which can vary in function of the neighborhood. 
If we say that it is {\em geometric} this includes the condition that it is globally an $n$-stack for a fixed $n$. 

The objects which arise here are naturally $(\infty ,1)$-stacks (which by convention p. \pageref{conven} we call just ``stacks''), 
that is they include morphisms which are not necessarily
invertible. For these, we make the somewhat {\em ad hoc} definition that a stack $\Mm $ is {\em locally geometric} if, on the one hand its $(\infty ,0)$-interior
(which uses only the invertible morphisms) 
is locally geometric in the above sense, and on the other hand for any two sections $E,F:S\rightarrow \Mm $, the $(\infty ,0)$-stacks $\Hom _{\Mm}(E,F)\rightarrow S$
are geometric. In the cases occuring in the present paper, these $\Hom$ stacks will be geometric because they come as Dold-Puppe of perfect complexes over $S$. 
Thus, the main problem is geometricity of the $(\infty ,0)$-interior. 

\subsection{Moduli stacks of perfect complexes on a formal category}
\mylabel{modstacks}

Suppose $f: (X,\Ff )\rightarrow Y$ is a morphism of smooth type from a formal category of smooth type to a scheme $Y$.   
If $Z\rightarrow Y$ is a morphism from an affine scheme,
the fiber product $f^{-1}(Z):= (X\times _YZ, \Ff \times _Y Z)$ is then a formal category whose associated $1$-stack is the similar fiber product
and corresponding to a sheaf of rings $\Lambda _{f^{-1}(Z)}$. 
Fix an affine open covering $\Uu$ of $X$, subordinate to an affine covering of $Y$. For any $Z\rightarrow Y$ the set of $U\times _YZ$ forms an affine covering
denoted $\Uu _Z$ of $X\times _YZ$. 

\begin{proposition}
\mylabel{allsame}
The following stacks on $\Aff /Y$ are all naturally equivalent: 
$$
\uHom ([X /\! /\Ff ],\Perf ):(Z/Y)\mapsto \Hom _{(\infty ,1)St/\Aff}([X /\! /\Ff ]\times _YZ, \Perf ),
$$
$$
\Mm ((X,\Ff )/Y, \Perf ):(Z/Y)\mapsto  L\cpx ^{\rm inj}_{\Oo -{\rm perf}}(CRIS(f^{-1}(Z)), \Oo ) ,
$$
$$
\Mm (\Lambda /Y, \Perf ):(Z/Y)\mapsto L\cpx ^{\rm inj}_{\Oo -{\rm perf}}(\Lambda _{f^{-1}(Z)}),
$$
$$
\Mm ^{\rm wk}(\Lambda /Y, \Perf ):(Z/Y)\mapsto \widetilde{DP} G_{\Uu _Z}([V\mapsto \wpx ^{\rm bplf}(\Lambda _{f^{-1}(Z)}(V)) ]) .
$$
Furthermore they all map to $\Perf (X/Y)$ and the equivalences respect this map. The weak version $\Mm ^{\rm wk}(\Lambda /Y, \Perf )$ is fibrant 
relative to $\Perf (X/Y)$. 
\end{proposition}
\begin{proof}
The first equivalence is by Theorem \ref{cpxonstack}.
For the remaining ones, note first of all that the maps of functoriality are obtained as in Lemma \ref{functoriality}.
For the equivalences, apply Lemma \ref{lambdastack} and the discussion of \S \ref{globwkcpx}. These also give the stack conditions
(for $\Mm ^{\rm wk}(\Lambda /Y, \Perf )$ the stack condition is a consequence of its objectwise
equivalence with the other versions).  The fibrancy condition for
$\Mm ^{\rm wk}(\Lambda /Y, \Perf )\rightarrow \Perf (X/Y)$ is meant in the sense of the HBKQ model structure of \cite{HirschowitzSimpson}
where we use Bergner's model structure for simplicial categories over each object \cite{Bergner}. Objectwise fibrancy implies fibrancy
because of the stack condition.
\end{proof}

The reader could have mostly skipped sections \S \ref{sec-crys}, \ref{sec-from}, \ref{sec-lambdamodel} if willing to forgo 
$\uHom ([X /\! /\Ff ],\Perf )$ or $\Mm ((X,\Ff )/Y, \Perf )$
and work directly with $\Mm (\Lambda /Y, \Perf )$. 

We could also consider the more general situation of a morphism of smooth type between formal categories of smooth type 
$f: (X,\Ff )\rightarrow (Y,\Hh )$. We can similarly define the moduli stack $\Mm ((X,\Ff )/(Y,\Hh ), \Perf )$ which can be viewed either
as a stack over $CRIS(Y,\Hh )$, or else a stack on $\Aff$ together with a map to $[Y/\! /\Hh ]$. 
For any $Z\rightarrow [Y/\! /\Hh ]$ we should first replace by an open covering on which the map lifts to $Z_i\rightarrow Y$,
look at perfect complexes over the formal category fibers $f^{-1} (Z_i)$ using any of the equivalent versions, then go back down to
an $(\infty ,1)$-category associated to $Z/Y$ using once again the descent machinery of the previous section. 
The pullback by $Y\rightarrow [Y/\! /\Hh ]$ is the moduli stack for the fiber product formal category over $Y$. 
This kind of construction yields,
for example, the Gauss-Manin connection on the family of $\Mm (X_{y,DR},\Perf )$ as $y\in Y$ varies when $f:X\rightarrow Y$ is a smooth morphism. 

If $f:(X,\Ff )\rightarrow (Y,\Hh )$ is a morphism of smooth type relative to another base scheme $S$ then we get a morphism of stacks
$$
\rr  f_{\ast}: \Mm ((X,\Ff )/S, \Perf ) \rightarrow \Mm ((Y,\Hh )/S, \Perf ),
$$
as well as $\rr  f^{\ast}$ going in the other direction. Again this could also be done relative to a third formal category $(S, \Gg )$.

\subsection{Maurer-Cartan stacks}
\mylabel{mcstack}

Suppose $Z$ is an augmented dga
over a commutative ring $R$.  An {\em MC element} is $\eta \in Z^1$ such that $d(\eta ) + \eta ^2=0$.
Suppose we are given an ideal $J\subset Z$ and an MC element $\eta _0\in Z/J$ which commutes with $1_{Z/J}$. 
Denote by $p:Z\rightarrow Z/J$ the projection. A {\em normalized MC element}
is an MC element $\eta$ such that $p(\eta )=\eta _0$. We define an augmented dg category of normalized MC elements
as follows. 

If $\eta$ and $\varphi$ are two normalized MC elements, the differential
$$
d_{\eta ,\varphi }(a):= d(a) + \varphi a - (-1)^{|a|}a\eta 
$$
has square zero. It preserves the subcomplex 
$$
J_+ := \ker (Z\rightarrow (Z/J)/\cc \cdot 1_{Z/J}) \hookrightarrow Z.
$$
Define the $\Hom$ complex by 
$$
\Hom _{\MC (Z,J,\eta _0)}(\eta , \varphi ):= \left( J_+  , d_{\eta ,\varphi }\right) .
$$
The projection $p: Z\rightarrow Z/J$ sends $J_+$ to
$\cc \cdot 1_{Z/J}$ and this provides an augmentation $\epsilon : \Hom _{\MC (Z,J,\eta _0)}(\eta , \varphi )\rightarrow \cc$. 
The composition of morphisms is given by the algebra structure on $Z$, and
it is compatible with the differentials. We obtain an augmented dg category $(\MC (Z,J,\eta _0), \epsilon )$.
Applying the affine Dold-Puppe construction (page \pageref{affineDP}) we get a simplicial category $\widetilde{DP}({\MC}(Z,J,\eta _0),\varepsilon )$.

\begin{definition} 
\mylabel{qnfilt}
Suppose $R$ is a $\cc$-algebra of finite type, and $Z$ is a differential graded algebra over $R$, with a decreasing filtration
$\{ J^kZ\}$ 
compatible with the differential and product. Say that $(Z,J)$ is {\em quasi-nilpotent} (resp. {\em nilpotent}) if 
there exists $k_0$ such that $J^kZ/J^{k+1}Z$ is acyclic for $k\geq k_0$, and furthermore $Z$ is complete with respect to
the filtration (resp. $J^k=0$). 
\end{definition}

This hypothesis provides a certain type of nilpotence which replaces the use of Artinian local rings as coefficients
in \cite{GoldmanMillson} for example. One consequence is homotopy invariance such as in \cite{GoldmanMillson2}.

\begin{lemma}
\mylabel{truncatek}
If $Z$ has a quasi-nilpotent filtration $J^{\cdot}$ and $k_0$ is chosen such that
$J^k/J^{k+1}$ is acyclic for $k\geq k_0$, then $Z/J^{k_0}$ has a nilpotent filtration.
If $\eta _0\in Z/J^1$ is an MC element, then the projection $p_{k_0}:Z\rightarrow Z/J^{k_0}$
induces a quasiequivalence of augmented dg categories
$$
\MC (Z,J^1,\eta _0)\stackrel{\sim}{\rightarrow} \MC (Z/J^{k_0}, J^1/J^{k_0},\eta _0).
$$
\end{lemma}
\begin{proof}
By the argument of Schlessinger-Stasheff-Deligne-Goldman-Millson
\cite{GoldmanMillson} \cite{GoldmanMillson2}, we can solve the necessary equations 
(either for proving quasi-fully-faithfulness or essential surjectivity) successively and inductively
to show the equivalence at the $k$-th stage for $k\geq k_0$
$$
\MC (Z/J^k,J^1/J^k,\eta _0)\stackrel{\sim}{\rightarrow} \MC (Z/J^{k_0}, J^1/J^{k_0},\eta _0).
$$
By the completeness hypothesis included in \ref{qnfilt} these solutions, chosen appropriately, fit together into a
solution for $Z$. 
\end{proof}

\begin{corollary}
\mylabel{invarianceMCstacks}
Suppose $\psi : Q\rightarrow Z$ is a morphism of differential graded $R$-algebras, such that both $Q$ and $Z$ have 
quasi-nilpotent filtrations
and $\psi$ preserves the filtrations. Suppose $\eta _0\in Q/J^1Q$ and $\varphi _0=\psi (\eta _0)\in Z/J^1Z$ are MC elements. 
If $\psi$ is a filtered quasiisomorphism except possibly at $J^0/J^1$, i.e. $Gr^m_J(\psi): Gr^m_J(Q)\stackrel{\sim}{\rightarrow} Gr^m_J(Z)$ is a quasiisomorphism for all $m\geq 1$,
then $\psi$ induces an equivalence
of augmented dg categories from ${\MC}(Q,J^1Q, \eta _0 )$ to ${\mathcal MC}(Z, J^1Z,\varphi _0 )$.
\end{corollary}
\begin{proof}
Choose $k_0$ sufficiently large for both $Q$ and $Z$, and note that $\psi$ induces an equivalence
$$
\MC (Q/J^{k_0}Q,J^1Q/J^{k_0}Q,\eta _0)\stackrel{\sim}{\rightarrow} \MC (Z/J^{k_0}Z, J^1Z/J^{k_0}Z,\varphi  _0)
$$
by the filtered quasiisomorphism hypothesis and the standard lifting procedure of \cite{GoldmanMillson} which terminates at $k_0$. 
Note that in the definition of $\MC$, the projections of everything modulo $J^1$ 
are constrained, both in the definition of MC element and in the definition of weak mapping complex. This is why we don't
need the quasiisomorphism  hypothesis for $J^0/J^1$. Applying Lemma \ref{truncatek} on both sides we get a commutative
square where three of the maps are quasiequivalences, so the fourth one is too.
\end{proof}

If $B$ is a commutative $R$-algebra we can make the $\MC$ construction for $Z\otimes _RB$,
then apply the affine $\widetilde{DP}$ construction. 
Let ${\calMC}(Z,J,\eta _0  )$ denote the stack on $\Aff / \Spec (R)$, associated to the prestack
$$
\Spec (B) \mapsto \widetilde{DP}({\MC}(Z\otimes _RB, J\otimes _RB, \eta _0\otimes 1_B),\varepsilon ).
$$
If $Z$ has a nilpotent filtration then so does $Z\otimes _RB$ (this isn't exactly true for a quasi-nilpotent grading;
we would instead have had to use a completed tensor product).

\begin{theorem}
\mylabel{finiteMCgeom}
Suppose that $Z$ has a nilpotent filtration, such that each $J^kZ/J^{k+1}Z$ is 
composed of free $R$-modules of finite rank, concentrated in cohomological degrees $\geq -n$. 
Suppose that $\eta _0\in Z/J^1$ is an MC element. 
Then ${\calMC}(Z,J^1,\eta _0 )$ is a geometric $n+1$-stack over $\Spec  (R)$.
\end{theorem} 
{\em Proof:}
The $Hom$ complexes of the Maurer-Cartan stack are obtained by affine Dold-Puppe of perfect complexes with amplitude in $[-n,\infty ]$.
Since $DP$ truncates at $0$, the unboundedness in positive degrees is not a problem. By \cite{geometricN} these
$Hom$ complexes are geometric. It remains to be seen how to define the smooth chart. The equation
$d(\eta ) + \eta ^2=0$ together with $\eta \mapsto \eta _0 \mod J^1$ defines a closed subvariety $V$ of an affine space of finite dimension over $R$. 
The tautological $\eta \in Z^1\otimes _R\Oo (V)$ is a map
$V\rightarrow  {\calMC}(Z,\varepsilon )$ which
is obviously surjective. To prove that it is smooth, suppose $g:Y\rightarrow {\calMC}(Z,\varepsilon )$
is another map from an affine scheme $Y$ corresponding to a Maurer-Cartan element $\varphi \in Z^1\otimes _R\Oo (Y)$.
We want to show that the map
\begin{equation}
\label{fiprodeq}
Y\times _{ {\calMC}(Z,\varepsilon )}V \rightarrow Y
\end{equation}
is smooth. The fiber product is geometric, by geometricity of the $Hom$ complexes in ${\calMC}$, 
so we just need to provide it with a chart smooth over $Y$.
The free $R$-module $Z^0$ corresponds to an affine space over $\Spec (R)$. Define $C$ to be the set of all
$(y, \alpha )\in Y\times _{\Spec (R)}Z^0$ with $\alpha \equiv 1 \mod J$. For any $(y, \alpha )$, there is a unique
solution $\psi \in V_{g(y)}$ with 
$$
d(\psi ) + \psi ^2=0
$$
\begin{equation}
\label{alphaeq}
d\alpha + (\alpha \psi - \varphi \alpha ) = 0.
\end{equation}
This gives a map $\Psi: C\rightarrow Y\times _{ {\calMC}(Z,\varepsilon )}V$. The 
scheme $C$ is clearly smooth over $Y$, so we just have to prove that $\Psi$ is a chart for the fiber product.

The two maps in the fiber product $Y\times _{ {\calMC}(Z,\varepsilon )}V$ 
correspond to $\varphi$ and $\eta$ respectively, and the fiber product is
the total space over $Y\times V$ of the affine $DP(\Hom _{\MC (Z,J,\eta _0)}(p_1^{\ast}\varphi , p_2^{\ast}\eta ),\varepsilon )$
consisting of morphisms projecting to $1$ by $\varepsilon$. Such morphisms are automatically invertible. 
This total space is a geometric stack (an affine Dold-Puppe of a perfect complex)  whose standard chart \cite{geometricN}
is exactly the set of solutions $(y,\psi ,\alpha )$ of the equations \eqref{alphaeq} above. Thus, $\Psi$ is isomorphic to the standard chart
for the fiber product.   
\eop

The $\MC$ construction commutes with globalization, taking into consideration the augmentations. 
If $U\mapsto Z_U$ is a presheaf of dga's on $\Uu$ with filtrations $J^{\cdot}Z_U$ then we obtain a filtration 
of the globalization defined by
$$
J^{k}G_{\Uu}[U\mapsto Z_U]:= G_{\Uu}[U\mapsto J^kZ_U].
$$
Globalization commutes with associated-graded, so it respects filtered quasiisomorphisms. 
Suppose $\eta _0$ is an MC element in $G_{\Uu}[U\mapsto Z_U] /J^1$. It projects to an MC element $\eta _{0,U}$ in
each $Z_U/J^1$. 

\begin{lemma}
\mylabel{MCglobcommute}
There is a natural equivalence of augmented dg categories
$$
\MC (G_{\Uu}[U\mapsto Z_U], J^1, \eta _0) \stackrel{\sim}{\rightarrow} 
G^{\rm eq}_{\Uu} \left( [U\mapsto \MC (Z_U, J^1Z_U, \eta _{0,U})] ,\epsilon  \right)  .
$$
\end{lemma}
\begin{proof}
Ojects on the left are collections $\eta (U_0,\ldots , U_k)$ indexed by strictly increasing sequences
$U_0\subset \cdots \subset U_k$ in $\Uu$ satisfying an MC equation of the form
$d\eta + \eta ^2=0$. Decompose into $\eta '$ which is the collection of $\eta (U_0)$, and $\eta ''$ which is
the collection of $\eta (U_0,\ldots , U_k)$ for $k\geq 1$. Then $\eta '$ can be viewed as a collection of objects
in the $\MC (Z_U, J^1Z_U, \eta _{0,U})$, while $\eta ''$ can be viewed as an MC element associated to this collection of
objects so as to define an element $(\eta ',\eta '')$ of the dgc globalization $G^{\rm eq}_{\Uu}$. The morphism complexes on 
both sides are the same, in this point of view. 
\end{proof}

\subsection{Explicit description of the homotopy fiber}
\mylabel{fibexplicit}

Recall that the simplicial category fiber of a fibration in dg-cat is the affine Dold-Puppe
of the augmented dg-fiber. We can combine the previous discussions to get an explicit identification. 

\begin{theorem}
\mylabel{mainident}
Suppose $E$ is a globally defined perfect complex expressed as a bounded complex of locally free $\Oo _X$-modules. 
There is a natural equivalence of simplicial categories between the homotopy fiber of 
$$
L\cpx ^{\rm inj}_{\Oo -{\rm perf}}(\Lambda )\rightarrow L\cpx ^{\rm inj}_{\Oo -{\rm perf}}(\Oo _X )
$$
over $E$, and the affine Dold-Puppe of the Maurer-Cartan category of the \v{C}ech globalization of the 
Hochschild-style dga
$$
\widetilde{DP}\left( \MC \left( G _{\Uu} [U\mapsto Q(E_U,E_U)], J^1, 1_{E}(1) \right) , \epsilon \right) .
$$
\end{theorem}
\begin{proof}
As discussed subsequent to Lemma \ref{globfib}, see also Lemma \ref{dpfiber},
the homotopy fiber of the map in question is equivalent to the Dold-Puppe of 
the augmented dg fiber of the map \eqref{mainmap} at $E$, which we shall denote simply by 
$\Fib ^+/E$ for brevity. 

In turn the augmented dg fiber  is
the augmented globalization of the functor of augmented dg-fibers (Lemma \ref{globfib})
$$
\Fib ^+/E =
G ^{\rm eq}_{\Uu}\left( [U\mapsto \Fib ^+( \wpx ^{\rm bplf}_0(\Lambda (U))\rightarrow DG\cpx ^{\rm bplf}(\Oo (U)))/E)],\epsilon \right) .
$$
For any $U$, the augmented dg fiber 
$$
\Fib ^+(( \wpx ^{\rm bplf}_0(\Lambda (U))\rightarrow DG\cpx ^{\rm bplf}(\Oo (U)))/E)
$$
is the relative Maurer-Cartan augmented dg category of the filtered dg algebra $Q(E_U,E_U)$ with respect to the
MC element $1_{E_U}(1) \in Q(E_U,E_U)/J^1$. 
We obtain the identification 
$$
\Fib ^+/E = G _{\Uu}^{\rm eq}\left( [U\mapsto \MC (Q(E_U,E_U); 1_E(1))] , \epsilon \right) .
$$
Commutation of the $\MC$ construction with globalization in Lemma \ref{MCglobcommute} gives
$$
\Fib ^+/E  = \MC \left( G _{\Uu} [U\mapsto Q(E_U,E_U)], J^1, 1_{E}(1) \right) .
$$
\end{proof}

\subsection{The Hochschild-Kostant-Rosenberg theorem}
\mylabel{sub-hkr}

If $F$ is a complex of strict $L$-modules, 
$\Hom _A(P,F)= \Hom _L(L\otimes _AP,F)$ for any complex of $A$-modules $P$. In particular, for a complex of $A$-modules $E$,
$$
Q(E,F) = \Hom _A(TL\otimes _AE,F) = \Hom_L(L\otimes _ATL \otimes _AE,F).
$$
Suppose $(E,\eta )$ is an MC object, and let $\phi$ be the MC element for $F$ corresponding to its strict
$L$-module structure (Lemma \ref{weakify}). 
Define a differential $d^+_{\eta}$ on $L\otimes _ATL\otimes _AE$ as follows:
$$
d^+_{\eta}(u_0\otimes  \cdots \otimes u_k\otimes e):=
\sum _{i=1}^{k-1}(-1)^{i+1} u_1\otimes 
\cdots \otimes (u_iu_{i+1})\otimes \cdots \otimes e 
$$
$$
+ \sum (-1)^m u_0\otimes \cdots \otimes u_m \otimes \eta (u_{m+1},\ldots , u_k;e).
$$
Thus $(L\otimes _ATL \otimes _AE, d^+_{\eta})$ is a complex of strict $L$-modules.

\begin{lemma}
\mylabel{identif}
Under the hypothesis that $F$ is a complex of strict $L$-modules the above construction identifies
$$
(Q(E,F),d_{\eta , \phi}) = \Hom ^{\cdot}_L((L\otimes _ATL \otimes _AE, d^+_{\eta}), F)
$$
where the right hand side is provided with the natural differential combining $d^+_{\eta}$ and $d_F$. 
\end{lemma}
\begin{proof}
The differentials are the same by construction.
\end{proof}

In order to apply Theorem \ref{finiteMCgeom} to prove geometricity, 
we need to show that what happens with complexes of $\Lambda$-modules only depends on a finite piece $J_k\Lambda$. 
This is a version of the Hochschild-Kostant-Rosenberg theorem. The proof is by reduction to the associated-graded, so we  need to understand what
happens for a polynomial ring. 

It is instructive to consider separately the case $L=A$. An element $q\in Q(E,F)$ may then be considered as a collection 
of $A$-linear functions $E\rightarrow F$,
$$
q^{m}(e)= q(1,\ldots , 1; e)\in F
$$
where there are $m$ $1$'s in the expression. The differential is
$(d_Qq)^m=q^{m-1}$ if $m$ is even, or $(d_Qq)^m=0$ if $m$ is odd.
The juxtaposition composition is given by
$(pq)^m = \sum _{i=0}^m p^iq^{m-i}$.
The projection 
$$
P_0: Q(E,F)\rightarrow \Hom ^{\cdot}_A(E,F), \;\;\; q\mapsto q^0
$$
is a surjection with an obvious splitting compatible with juxtaposition and the differential, and $\ker (P_0)$ is acyclic. 

Next, consider the case where $L=A[x_1,\ldots , x_m]$ is a commutative polynomial ring over $A$.
Then we have an augmentation $v:L\rightarrow A$ (evaluation at $x_i=0$) so any complex of $A$-modules can be considered as
a complex of $L$-modules with $L$ acting through $v$. In particular, both $E$ and $F$ have structures of strict $L$-modules,
with their associated MC elements which we denote $\varphi$ and $\phi$
respectively. We can apply the identification of Lemma \ref{identif}:
$$
(Q(E,F),d_{\varphi , \phi })= \Hom _L((L\otimes _ATL\otimes _AE , d^+_{\varphi}),F).
$$
We have a quasiisomorphism of complexes of $L$-modules
$$
(L\otimes _ATL\otimes _AE , d^+_{\varphi}) \stackrel{\sim}{\rightarrow} (E,d_E),
$$
shown as in the proof of Lemma \ref{wkLacyclic} by combining into a complex of the form
$TL\otimes _AE$ which is acyclic using a homotopy obtained by putting in $1$'s. 
The complex $Q(E,F)$ as identified above therefore calculates the groups ${\rm Ext}^i_L(E,F)$. 

\begin{lemma}[\cite{HKR}]
\mylabel{polycase}
Suppose $L=A[x_1,\ldots , x_m]$ is a commutative polynomial ring over $A$. Suppose
$E$ and $F$ are fixed perfect complexes of $A$-modules, realized as bounded complexes of locally free $A$-modules.
Provide $Q(E,F)$ with the differential corresponding to the structures of $L$-modules on $E$ and $F$ via the evaluation at $x_i=0$.
For any $k\geq m+1$ the complexes $J^kQ(E,F)/J^{k+1}$ and $J^kQ(E,F)$
are acyclic. 
\end{lemma}
\begin{proof}
The direct sum decomposition of $TL$ according to total degree of polynomials, gives a direct product
decomposition of $\Hom ^{\cdot}(L\otimes _ATL\otimes _AE,F)$. This in turn leads to a direct product decomposition of 
${\rm Ext}^i_L(E,F)$. The pieces may be identified by looking at 
the scalar action of multiplying the coordinates $x_1,\ldots , x_m$ by $t\in \cc ^{\ast}$.
The $k$-th graded piece in this decomposition is the cohomology of the complex $J^kQ(E,F)/J^{k+1}$. 

On the other hand, we can use a standard Koszul resolution to calculate ${\rm Ext}^i_L(E,F)$, and the pieces of the 
Koszul resolution have homogeneity of degree $\leq m$, the biggest homogeneity being from the last factor $\bigwedge ^m _L(L^m)$.
This shows, for one thing, that the action of $t\in \cc ^{\ast}$ is algebraic, which shows that there can only be finitely
many nontrivial terms in the direct product. Then consideration of the degrees of homogeneity shows that only the terms of
degree $\leq m$ contribute. We conclude that for $k\geq m+1$ the complex $J^kQ(E,F)/J^{k+1}$ has no cohomology. 
The complex $J^kQ(E,F)$ is a direct product of the subquotient complexes in degrees $\geq k$ so it is also acyclic whenever $k\geq m+1$. 
\end{proof}

\subsection{Acyclicity at high orders}
\mylabel{sub-acyclic}

Recall that $L$ is almost-polynomial, in other words $Gr^J(L)$ is a polynomial ring, due to the fact that
$L$ comes from a formal category of smooth type.

\begin{proposition}
\mylabel{acyclic}
Let $m$ be the rank of the locally free $A$-module
$J_1L/J_0L$. Let $(E,\eta )$ and $(F,\xi )$ be MC objects. For any $k\geq m+1$ the complexes
$(J^kQ(E,F)/J^{k+1}, d_{\eta , \xi })$ and $(J^kQ(E,F), d_{\eta , \xi })$ are acyclic.  
\end{proposition}
{\em Proof:} 
Recall that Wodzicki \cite{Wodzicki} gave a spectral sequence argument for this type of result,
reducing it to the case of a polynomial ring \cite{HKR}.  
The Hochschild-style complex $TL$ we are using here differs from the
one in \cite{Wodzicki} because we are contracting inner $A$-module structures where as the complex considered in
\cite{Wodzicki} involves tensoring over the ground field which in our case is $\cc$. Nonetheless,
the basic principle of the proof remains the same: we can take the associated-graded for our filtration,
and then the calculation reduces to that of a polynomial algebra over $A$ where we can apply
\cite{HKR}. In our current notations this argument comes down to noting that 
the complex $(J^kQ(E,F)/J^{k+1}, d_{\eta , \xi })$ is the same as the corresponding one for the polynomial ring
$Gr^J(L)$, then apply Lemma \ref{polycase}. Once we know that the elementary subquotients are acyclic it follows that the
bigger subquotients of the form $(J^kQ(E,F)/J^{k+p}, d_{\eta , \xi })$ are acyclic by exact sequences. Then
$J^kQ(E,F) = \lim_{\leftarrow , p}J^kQ(E,F)/J^{k+p}$ and the transition maps in the inverse system are surjective,
so $(J^kQ(E,F), d_{\eta , \xi })$ is acyclic. 
\eop

\begin{corollary}
\mylabel{quasiff}
If $m:={\rm rk}\, \Omega ^1_{(X,\Ff )}$, for any $k\geq m+1$ the morphism
$$
( Q(E,F),d_{\eta , \xi} ) \rightarrow ( Q(E,F)/J^k,d_{\eta , \xi} )
$$
is a filtered quasiisomorphism. The filtered dga $Q(E,E)$ is quasi-nilpotent \ref{qnfilt},
and the fiber in Theorem  \ref{mainident}
is equivalent to 
$$
\widetilde{DP}\left( \MC \left( G _{\Uu} [U\mapsto Q(E_U,E_U)/J^k], J^1, 1_{E}(1) \right) , \epsilon \right) .
$$
\end{corollary}
\eop

Although we don't need it, one can define a category like $\wpx ^{\rm bplf}_0(L)$ but using the quotients $Q(E,F)/J^k$ both for
the MC elements $\eta$ and for the weak mapping complexes. The previous corollary says that the functor from
$\wpx ^{\rm bplf}_0(L)$ to this quotient category is a quasiequivalence of dg categories
relative to and fibrant over $DG\cpx ^{\rm pblf}(A)$.

\subsection{A finite-dimensional replacement}
\mylabel{finiterep}

The last step of the proof is to reduce from a \v{C}ech complex to a complex involving sections with a bounded
number of poles along the complementary divisors of the affine open sets. This will give a finite dimensional complex. 

Suppose $X\rightarrow \Spec  (R)$ is a smooth projective map
with $R$ a commutative $\cc$-algebra of finite type.

Suppose $E$ is a perfect complex  which is bounded and 
whose components are locally free $\Oo _X$-modules of finite rank. Note that any perfect complex
has a representative of this form \cite{ThomasonTrobaugh}, SGA VI. 

In the global situation we can look at $\uQ (E,E)= \uHom^{\cdot} _{\Oo _X}(T^{\cdot}\Lambda \otimes _{\Oo _X}E,E)$
which is an inverse system of complexes of coherent sheaves on $X$. It has a filtration by $J^k\uQ (E,E)$ (corresponding
to the filtration $J_p\Lambda$ as before) and the quotients $\uQ (E,E)/J^k$ are locally free $\Oo _X$-modules of finite rank, in each
cohomological degree. It is bounded below but not above. 

On an affine open set $U$, we have $\uQ (E,E)(U)= Q(E_U,E_U)$ as considered previously. 
For any $k\geq m+1$, the map $\uQ (E,E)(U)\rightarrow \uQ (E,E)(U)/J^k$ induces a
quasiequivalence of dg categories of MC elements by \ref{quasiff} and \ref{invarianceMCstacks}. 

Recall from \S \ref{sheafglob} the {\em sheaf-theoretic globalization} 
$\Gg _{\Uu}(\uQ (E,E)/J^k)$ which is a sheaf of dga's quasiisomorphic to $\uQ (E,E)$, the global sections of which 
are identical to $G_{\Uu}[U\mapsto Q(E_U, E_U)/J^k]$. However,
its terms are direct images of coherent sheaves from open sets into $U$, so they
are not coherent sheaves over $X$. We would like to replace this by a coherent sheaf of dga's.

Fix an increasing positive function $\underline{m}(k)$ of $k=0,1,\ldots $. For any divisor $D$ on $X$ relative to $Spec(R)$, define
$$
\uQ (E,E)(\underline{m}^JD)/J^k := \sum _{p=0}^{k-1}
(J^p\uQ (E,E)/J^k)\otimes _{\Oo _X}\Oo _X(\underline{m}(p)D).
$$
The sum can be taken as a subcomplex of $(\uQ (E,E)/J^k)\otimes _{\Oo _X}\Oo _X(\underline{m}(k))$. 
To put it another way, the sections here are those which have poles of order $\underline{m}(p)$ when projected into 
$\uQ (E,E)/J^p$. 
Denote by a non-script letter the complex of global sections:
$$
Q (E,E)(\underline{m}^JD)/J^k := \Gamma (X,\uQ (E,E)(\underline{m}^JD)/J^k ).
$$
Notice that $Q (E,E)(\underline{m}^JD)/J^k$ is a subcomplex of $\uQ (E,E)(U)/J^k$ where $U=X-D$. 

Assume that $\underline{m}(p)\underline{m}(q)\leq \underline{m}(p+q)$. In particular this requires $\underline{m}(0)=0$.
With this assumption, $\uQ (E,E)(\underline{m}^JD)/J^k $ remains a sheaf of algebras on $X$ and 
$Q (E,E)(\underline{m}^JD)/J^k$ is a dga. 

This construction can be mixed in with the \v{C}ech globalization. Suppose that our open covering 
$\Uu$ consists of open sets  $U_i=X-D_i$ where $D_i$ is the divisor of a section $s_i$ of a fixed very ample
line bundle $\Oo _X(1)$. Thus the multiple intersections $U_I=U_{i_1,\ldots , i_m}$ are the complements
of the divisors $D_I$ given by sections $s_I$ of $\Oo _X(m)$ where $m=|I|$.
Let ${\bf D}$ denote this collection of divisors.

In the definition of the sheaf-theoretical globalization $\Gg _{\Uu}(\uQ (E,E)/J^k)$, 
replace $j_{U_I/X, \ast}\uQ (E_{U_I}, E_{U_I})/J^k$ by 
an object $\uQ (E,E)(\underline{m}^J{\bf D} )/J^k$ defined by making the above replacement for each $U_I$. Call the resulting complex of sheaves
$\Gg _{\underline{m}{\bf D}}(\uQ (E,E)/J^k)$, and its complex of global sections is denoted
$G_{\underline{m}{\bf D}}(\uQ (E,E)/J^k)$. It
fits in between $\uQ (E,E)/J^k$ and the full \v{C}ech globalization 
$$
\uQ (E,E) /J^k\rightarrow \Gg _{\underline{m}{\bf D}}(\uQ /J^k)\rightarrow \Gg _{\Uu}(\uQ (E,E)/J^k).
$$
These are quasiisomorphisms of complexes of sheaves (this is best understood by looking at the example in \S \ref{sheafglob}). 
We can assume that $\underline{m}(k)$ is big enough so that the terms of $J^1\Gg _{\underline{m}{\bf D}}(Q)$
are acyclic and their global sections are locally free over the base scheme. 
Then the local quasiisomorphism on the right induces a quasiisomorphism on 
the dga's of global sections.

\begin{proposition}
\mylabel{finitereplacement}
Consider a bounded family of
strictly perfect complexes of $\Oo _X$-modules $E$. We can fix $\underline{m}$ uniformly to get 
a differential graded $R$-algebra $G _{\underline{m}{\bf D} /X}(\uQ (E,E)/J^k)$ which is locally free of finite rank over $R$,
such that the map 
$$
G _{\underline{m}^J{\bf D} }(\uQ (E,E)/J^k)\rightarrow  G _{\Uu }[U\mapsto Q (E_U,E_U)/J^k]
$$
induces quasiisomorphisms on all graded pieces $J^m/J^{m+1}$ except $J^0/J^1$, and the same holds after any change of the base ring $R$.
\end{proposition}
\begin{proof}
It suffices to choose $\underline{m}$ so that the global sections of the graded pieces are locally free $R$-modules, compatible
with base change of the ring $R$, and the higher direct images to $Spec(R)$ vanish. This should hold starting with $J^1/J^2$. 
Then the map 
$\Gg _{\underline{m}{\bf D}}(\uQ (E,E)/J^k)\rightarrow \Gg _{\Uu}(\uQ (E,E)/J^k)$ induces a filtered quasiisomorphism
on global sections, again after any change of base and except at $J^0/J^1$. 
\end{proof}

\subsection{The geometricity theorem}
\mylabel{geomtheorem}

\begin{theorem}
\mylabel{thethm}
Suppose $(X,\Ff )$ is a formal category of smooth type over a base scheme $Y$ such that the underlying
map of schemes $X\rightarrow Y$ is smooth. Let $\Lambda$ be the corresponding
sheaf of rings of differential operators. 
The moduli stack \ref{allsame} $\Mm ((X,\Ff )/Y, \Perf )\approx \Mm (\Lambda /Y, \Perf )$ is a locally geometric stack over $Y$. 
\end{theorem}
\begin{proof}
The $Hom$ stacks are Dold-Puppe of truncated $Hom$ complexes between objects. The
$Hom$ complexes are perfect as can be seen using the standard argument with a de Rham-Koszul resolution. 
Thus it suffices to consider the $(\infty ,0)$-interior. Given a morphism of 
$(\infty ,0)$-stacks such that the target is geometric and the pullback to any scheme mapping to the
target is also geometric, then the source is geometric \cite[Corollary 2.6]{geometricN}. Apply this to the forgetful functor
$\Mm (\Lambda /Y, \Perf )\rightarrow \Perf (X/Y)$ using \cite{ToenVaquie} which we have stated as Theorem \ref{tv} above. Given a scheme mapping to $\Perf (X/Y)$ we can
pull back to there and call that new sheme $Y$, that is to say we can assume that we are given a perfect 
complex $E$ on $X$ and we have to show geometricity of the fiber 
$$
\Mm (\Lambda /Y; E):= \Mm (\Lambda /Y, \Perf )\times _{\Perf (X/Y)}\{ E \} .
$$
This is the stack which associates to an affine scheme $Z\rightarrow Y$ the homotopy
fiber of 
$$
L{\cpx}^{\rm inj}_{\Oo -{\rm perf}}(\Lambda |_{X_Z}) \rightarrow \Perf (X_Z)
$$
over $E_Z:=E|_{X_Z}$ where $X_Z:= X\times _YZ$. 

Fix an affine open covering $\Uu$ of $X$. We can be assuming that $Z$ is affine, so we get an affine open
covering $\Uu _Z$ of $X_Z$. By Theorem \ref{mainident} and Corollary \ref{quasiff}, 
$\Mm (\Lambda /S; E)(Z\rightarrow Y)$ is equivalent to 
$$
\widetilde{DP}\left( \MC \left( G _{\Uu _Z} [U\mapsto Q(E_U,E_U)/J^k], J^1, 1_{E_Z}(1) \right) , \epsilon \right) .
$$

We can assume that our open covering is given by a system of divisors ${\bf D}$. Note that by localizing over the base
and possibly throwing out some elements of the covering, we may assume that our divisors remain divisors after base change.
Thus, on any $X\times _YZ$ we get a system of divisors ${\bf D} (Z)$. The $E_Z$ vary in a bounded family so 
by Proposition \ref{finitereplacement} for an appropriate $\underline{m}$ uniform in $Z$, we have a filtered quasiisomorphism of dga's
(outside of $J^0/J^1$) 
$$
G _{\underline{m}^J{\bf D} (Z)}(\uQ (E_Z,E_Z)/J^k)\stackrel{\sim}{\rightarrow}  G _{\Uu _Z }[U\mapsto Q (E_U,E_U)/J^k] .
$$
Extending Corollary \ref{invarianceMCstacks} to stacks by applying it objectwise, we get that
$\Mm (\Lambda / Y; E)(Z\rightarrow Z)$ is equivalent to 
$$ 
\widetilde{DP}\left( {\MC}\left(  G _{\underline{m}^J{\bf D} (Z)}(\uQ (E_Z,E_Z)/J^k), J^1,1_{E_Z}(1) \right) ,
\varepsilon \right) . 
$$
The sheaf of dga's $\uQ (E_Z,E_Z)$ is the pullback to $X_Z$ of $\uQ (E,E)$ on $X$, and
$G _{\underline{m}^J{\bf D} (Z)}(\uQ (E_Z,E_Z)/J^k)$ is the pullback to $Z=\Spec (B)$ of 
the dga $G _{\underline{m}^J{\bf D} }(\uQ (E,E)/J^k)$ over $Y=\Spec (R)$. The filtration here is nilpotent. 
We get an identification 
$$
M(\Lambda / Y; E)\cong {\calMC}\left(  G _{\underline{m}^J{\bf D} }(\uQ (E,E)/J^k), J^1,1_E(1) \right) .
$$
Proposition \ref{finitereplacement} says that $G _{\underline{m}^J{\bf D} }(\uQ (E,E)/J^k)$ satisfies the hypothesis of Theorem \ref{finiteMCgeom}, therefore 
$\Mm (\Lambda / Y; E)$
is geometric. 
\end{proof}

An obvious question for further study is to extend this to derived stacks \cite{Kontsevich}
\cite{KapranovD} \cite{CFK} \cite{hag2} \cite{LurieD}.


\section{Application to the Hodge filtration}
\mylabel{rees-sec}

Suppose $X$ is a smooth projective variety over $\Spec  (\cc )$. The de Rham formal category $\Ff _{DR}\subset X\times X$
is the formal completion of the diagonal. Deformation to the normal cone \cite{Fulton} yields a deformation from $\Ff _{DR}$ to
the formal completion $\Ff _{Dol}$ of the zero-section in the tangent bundle of $X$ (whose structure of formal category is given by addition in the tangent bundle). 
This gives a $\Gm$-equivariant formal category
$X_{\rm Hod}= (X\times \aaa ^1 , \Ff _{\rm Hod})$ such that the morphism $\lambda : X_{\rm Hod}\rightarrow \aaa ^1$ is of smooth type \cite{hodgefilt}. The fiber
over any $\lambda \neq 0$ is isomorphic $X_{DR}=(X,\Ff _{DR})$, while the fiber over $\lambda = 0$ is $X_{Dol} = (X,\Ff _{Dol})$. 
The sheaf of rings of differential operators $\Lambda _{Hod}$ corresponding to $\Ff _{Hod}$ on $X\times \aaa ^1$ may be obtained by the Rees construction
as applied to $\Lambda _{DR}= \Dd _X$. 

By Theorem \ref{thethm}, $\Mm (X_{\rm Hod},\Perf ) := \uHom (X_{\rm Hod}/\aaa ^1 , \Perf )$ is a locally geometric stack. A section $S\rightarrow
\Mm _{\rm Hod}(X,\Perf )$ lying over a map $\varphi : S\rightarrow \aaa ^1$ is a complex $E$ of 
$\varphi ^{\ast}\Lambda _{\rm Hod}$-modules on $X\times S$, which is perfect over $\Oo _{X\times S}$. 
For any closed point $s\in S$ with $\varphi (s)\neq 0$ the restriction 
$E_s := E|_{X\times \{ s\}}$ is an $\Oo _X$-perfect complex of $\Lambda _{DR}= \Dd _X$-modules. Thus $E_s$ is a finite complex whose cohomology objects $H^i(E_s)$
are $\Dd _X$-modules locally free over $\Oo _X$, i.e. vector bundles with connection. 

For a closed point $s\in S$ with $\varphi (s)=0$, the restriction $E_s$ is an $\Oo _X$-perfect complex of $\Lambda _{\rm Dol}$-modules. Its cohomology objects are Higgs sheaves
which need not, however, satisfy any good conditions such as semistability, vanishing Chern classes, local freeness or even torsion-freeness. In order to obtain an object which should be reasonably considered as the {\em Hodge filtration} for $ \Mm _{DR}(X,\Perf )$ we should impose these conditions over points where $\varphi (s)=0$.

Define $\Mm _{Hod}(X,\Perf )$ to be the full substack of $ \Mm (X_{\rm Hod},\Perf )$ consisting of those sections $\phi : S\rightarrow  \Mm (X_{\rm Hod},\Perf )$
such that $(\ast )$ for any point $s\in S$ with $\lambda \circ \phi (s)=0$, the corresponding $\Oo _X$-perfect complex $E_s$ of $\Lambda _{\rm Dol}$-modules has cohomology 
objects $H^i(E_s)$ which are locally free Higgs bundles, semistable, with vanishing Chern classes on $X$. 

Fix a vector of positive integers $\underline{b} = \{ \underline{b}_i\} _{i\in \zz }$ such that $\underline{b}_i=0$ for all but a finite number of $i$. 
Recall that $\Perf (\underline{b})$ denotes the open substack of $\Perf$ consisting of perfect complexes whose Betti numbers are bounded by $\underline{b}_i$.
Let $ \Mm _{Hod}(X,\Perf (\underline{b}))$ be full substack of $ \Mm (X_{\rm Hod},\Perf (\underline{b}))$
defined by condition $(\ast )$ of the previous paragraph.

\begin{theorem}
\mylabel{hfthm}
The condition $(\ast )$ is an open condition in $ \Mm (X_{\rm Hod},\Perf )$, so it defines a locally geometric open substack.
It is covered by geometric stacks  $ \Mm _{Hod}(X,\Perf (\underline{b}))$. 
\end{theorem}
\begin{proof}
In order to show openness of the condition $(\ast )$, which only concerns $\lambda =0$, 
it suffices to consider the case when the base $S$ is a smooth affine curve, 
and when the projection $\varphi : S\rightarrow \aaa ^1$ is identically zero. In this case $\varphi ^{\ast}\Lambda _{\rm Hod} = 
\Lambda _{\rm Dol} (X\times S/S)$. 
Suppose $E$ is an $\Oo$-perfect complex of $\Lambda _{\rm Dol}$-modules on $X\times S /S$, satisfying condition $(\ast )$ at 
a point $s\in E$. Let $n$ be the upper limit for the interval
of amplitude of $E_s$, so $H^n(E_s)\neq 0$ but $H^i(E_s)=0$ for $i>n$. Possibly after replacing $S$ by a smaller neighborhood of $s$
we may assume that $n$ is the limit of the interval of amplitude for $E$, and also that the subobject of $\Oo _S$-torsion $H^n(E)_{\rm tors}$ 
lies over the point $s$. Let $t$ be a  local coordinate for $S$ at the point $s$, so $t$ acts nilpotently on $H^n(E)_{\rm tors}$.
The object $E_s = E\otimes _{\Oo _S} \cc $ may be calculated using the resolution 
$$
0\rightarrow \Oo _S \stackrel{t}{\rightarrow} \Oo _S \rightarrow \cc \rightarrow 0.
$$
We get a quasiisomorphism between $E_s$ and the total complex of $E\stackrel{t}{\rightarrow}E$. The spectral sequence for this double complex
gives an exact sequence
$$
H^{n-1}(E_s)\rightarrow H^n(E)_{\rm tors}\stackrel{t}{\rightarrow}H^n(E)_{\rm tors} \rightarrow H^n(E_s)\rightarrow \left( H^n(E)/{\rm tors}\right) _s \rightarrow 0.
$$
By the assumption of condition $(\ast )$ at the point $s$, the terms $H^{n-1}(E_s)$ and $H^n(E_s)$ are semistable Higgs bundles with vanishing Chern classes on 
$X$. 
An argument using the monodromy weight filtration for the nilpotent endomorphism $t$ implies that $H^n(E)_{\rm tors}$ is itself a semistable Higgs bundle with
vanishing Chern classes, and similarly for any of the subobjects of $t^k$-torsion, as well as 
the fiber of the torsion-free part $\left( H^n(E)/{\rm tors}\right) _s$. Note that $\left( H^n(E)/{\rm tors}\right)$ is a flat family of Higgs sheaves 
on $X\times S/S$, so this condition at the point $s$ implies the same at any other point $y$. 

If $V$ is a semistable Higgs bundle with vanishing Chern classes on 
$X\times S/S$ provided with a morphism $V\rightarrow H^n(E)$, we can lift it to a morphism $\psi : V [-n]\rightarrow E$ and look at $E':= {\rm Cone}(\psi )$.
As long as the interval of amplitude of $E$ has length at least $1$, the interval of amplitude of $E'$ will be no bigger. Using appropriate choices of $V$ and $\psi$,
made possible by the above observations, we can successively reduce the length of the torsion in $H^n(E)$ and then take out the torsion-free part. 

Proceeding by induction we obtain the statement that $E$ is obtained by successive mapping cones on semistable Higgs bundles with vanishing Chern classes. 
The maps involved should be calculated in the homotopically correct way, for example using injective resolutions. 
This devissage shows that condition $(\ast )$ holds at all points of $S$. We admitted replacement of $S$ by an open neighborhood of $s$ at the start, so this shows openness
of the condition $(\ast )$. It defines a locally geometric open substack $\Mm _{\rm Hod}(X,\Perf )\subset \Mm (X_{\rm Hod},\Perf )$. 

Making the convention that $\underline{b}_i=0$ for all but a finite number of $i$,
the open substack $ \Mm _{Hod}(X,\Perf (\underline{b}))$ is geometric. This is because the family of possible $E_s$ with Betti numbers less than $\underline{b}_i$
is bounded. 
\end{proof}

It is interesting to note that the above argument gives a structure result for sections of $ \Mm (X_{\rm Hod},\Perf )$ over a smooth curve. 

\begin{proposition}
If $S$ is a smooth curve, then etale-locally on $S$ any morphism $S\rightarrow  \Mm _{\rm Hod}(X,\Perf )$ projecting to $\varphi : S\rightarrow \aaa ^1$ 
comes from a double complex $E$ on $X\times S$ whose rows are injective resolutions of 
are vector bundles with semistable $\varphi $-connections on $X\times S/S$, with vanishing Chern classes
on the fibers. 
\end{proposition}
\begin{proof}
Use the same argument as before. The maps in the mapping cones can be realized as maps between injective resolutions and this leads to a double complex.
The function $\varphi$ may no longer be identically zero. In order to choose the objects $V$,
use the following observation: if $s\in S$ with $\varphi (s)=0$ and if $V_s$ is a semistable Higgs bundle with vanishing Chern classes on $X$,
then there exists a family $V$ on $X\times S/S$ of vector bundles with $\varphi$-connections extending $V_s$. Deligne's preferred sections of the
twistor space provide an analytic solution, and Artin approximation assures the existence of an algebraic solution after etale localization. 
\end{proof}

\begin{theorem}
Suppose $f: X\rightarrow Y$ is a smooth morphism between smooth projective varieties. Then the higher direct image defines a morphism
of geometric stacks
$$
\rr f_{\ast}: \Mm (X_{\rm Hod},\Perf )  \rightarrow \Mm (Y_{\rm Hod},\Perf ),
$$
whose fiber over $\lambda = 1$ is the higher direct image from complexes of $\Dd _X$-modules to complexes of $\Dd _Y$-modules. 
\end{theorem}
\begin{proof}
Use Lemma \ref{fcdirim}, and note that the condition $(\ast )$ is preserved by higher direct image \cite{families}. 
\end{proof}

The Deligne glueing operation works here as in the classical case \cite{hodgefilt}: we can glue $\Mm _{\rm Hod}(X,\Perf (\underline{b}))^{\rm an}$
to $ \Mm _{\rm Hod}(\overline{X},\Perf (\underline{b}))^{\rm an}$ to obtain the {\em Deligne-Hitchin twistor family} 
$$
 \Mm _{\rm DH}(X,\Perf (\underline{b}))\rightarrow \pp ^1.
$$
In order to make this precise we need a theory of analytic geometric stacks and the Riemann-Hilbert correspondence
$$
\Mm (X_{DR},\Perf (\underline{b}))^{\rm an} \cong \Mm (X_B,\Perf (\underline{b}))^{\rm an}.
$$
It would be good to investigate the notion of preferred section of the twistor space. The local deformation theory along such sections
should be governed by some kind of mixed Hodge complex \cite{hodge3} \cite{Hain2} \cite{namhs} \cite{MSaito0504}.


\end{document}